\documentclass[msom,sglanonrev]{informs4}

\OneAndAHalfSpacedXI

\usepackage{booktabs}           
\usepackage{multirow}           
\usepackage{longtable}          
\usepackage{diagbox}            

\usepackage{subcaption}         
\usepackage{caption}            
\usepackage{rotating}           

\usepackage{algorithm,algpseudocode}

\usepackage{siunitx}            
\usepackage{textcomp}           

\usepackage{tikz}
\usepackage{pgfplots}
\usepackage{pgfplotstable}
\pgfplotsset{compat=1.18}
\usetikzlibrary{shapes,arrows,positioning,backgrounds,calc}
\usetikzlibrary{arrows,decorations}

\DeclareMathSymbol{@}{\mathord}{letters}{"3B}

\definecolor{blues1}{RGB}{198, 219, 239}
\definecolor{blues2}{RGB}{158, 202, 225}
\definecolor{blues3}{RGB}{107, 174, 214}
\definecolor{blues4}{RGB}{49, 130, 189}
\definecolor{blues5}{RGB}{8, 81, 156}
\definecolor{blues6}{RGB}{2, 34, 78}

\pgfplotscreateplotcyclelist{mylist}{
{blues1},
{blues2},
{blues3},
{blues4},
{blues5},
{blues6},
}

\newcommand{\appendix}{\par\setcounter{section}{0}\gdef\thesection{Appendix \Alph{section}}}

\TheoremsNumberedThrough     
\ECRepeatTheorems            

\EquationsNumberedThrough    

\begin{document}

\TITLE{When to Match: A Cost-Balancing Principle for Dynamic Markets}

\ARTICLEAUTHORS{%
\AUTHOR{Jie Liu}
\AFF{Department of Industrial Engineering \& Decision Analytics, The Hong Kong University of
Science and Technology, Clear Water Bay, Hong Kong S.A.R., China,
\EMAIL{iejialiu@ust.hk}}
\AUTHOR{Hailun Zhang}
\AFF{School of Economics and Business Administration, Chongqing University, Chongqing, China,  \EMAIL{zhanghailun@cqu.edu.cn}}
\AUTHOR{Jiheng Zhang}
\AFF{Department of Industrial Engineering \& Decision Analytics, The Hong Kong University of
Science and Technology, Clear Water Bay, Hong Kong S.A.R., China,
\EMAIL{jiheng@ust.hk}}
}

\ABSTRACT{%
\textbf{Problem definition:}
Platforms in ridesharing, food delivery, and online gaming must decide not only whom to match but when: immediate matching cuts waiting, while delay thickens the market and improves match quality.
Because demand is hard to forecast, the right waiting window shifts continuously.
Fixed-window industry rules are simple but fragile, while forecast-based optimization models are brittle when assumptions fail.
This paper develops a matching rule that is as simple as industry practice yet carries a guarantee requiring no forecasts.
\textbf{Methodology/results:}
We study a model in which agents of several types are matched in groups drawing one agent from each type, waiting is costly, and matching costs fall as queues grow.
We propose the Cost-Balancing (CB) rule: match as soon as the waiting cost accumulated since the last match reaches a calibrated proportion of the current matching cost.
This principle is motivated by a steady-state fluid analysis of the matching system.
On any finite arrival stream delivering equal numbers of each type, CB calibrated for the worst case incurs at most twice the cost of an optimal clairvoyant policy that knows all future arrivals.
No deterministic online rule can guarantee a smaller factor, so CB is worst-case optimal, while greedy and fixed-threshold policies can perform arbitrarily worse than this benchmark.
The guarantee extends to matches with fixed heterogeneous consumption requirements.
In a game-matching experiment, CB reduces total cost by 3--8\% versus the industry-standard heuristic; in a food-delivery experiment, it reduces average delay by 14.5\% versus the best fixed-rule benchmark.
\textbf{Managerial implications:}
Platforms can manage match timing with Cost-Balancing, a simple, efficient, and robust rule that keeps waiting and matching costs in balance.
Its worst-case guarantee provides a safety net even in volatile conditions where fixed rules break down.
Responding to realized costs, the rule matches faster during surges and waits longer during lulls, without forecasts or retuning.
}


\maketitle

\section{Introduction}
\label{sec:intro}

Dynamic matching markets have become central to the modern service economy.
From transportation and logistics to entertainment and healthcare, the efficient allocation of resources is no longer a static planning problem but a continuous, real-time control challenge.
Platforms like Uber, Lyft, DoorDash, and competitive gaming servers operate in environments where supply and demand arrive stochastically, often with significant volatility.
In these fast-paced markets, the platform's central algorithm acts as the market maker, bearing the responsibility of deciding not just who matches with whom, but when these matches occur \citep{uber2023}.

This temporal dimension introduces a fundamental operational trade-off.
On the one hand, the need for responsiveness drives a preference for immediate matching.
Delays lead to explicit waiting costs: idle drivers burn fuel, hungry customers grow impatient, and patients suffer.
Excessive waiting also increases the risk of abandonment, leading to lost revenue and market contraction.
On the other hand, patience yields the benefit of market thickness.
By strategically delaying decisions, a platform can aggregate a critical mass of agents, transforming a sparse, inefficient market into a dense, efficient one where economies of scale can be exploited.
The magnitude of this trade-off varies across contexts, but its presence is universal in dynamic matching systems.

The business implications are substantial.
In on-demand delivery, the difference between profitability and loss often hinges on the batch rate, i.e., the average number of orders a courier can deliver per trip.
Immediate dispatch forces a one-to-one ratio, maximizing speed but crippling efficiency; batching can achieve higher ratios and reduce unit costs, but risks violating delivery guarantees and customer satisfaction.
In the multi-billion-dollar e-sports industry, player retention is driven by the flow state achieved in well-matched games.
An impatient algorithm that prioritizes short queues over skill parity creates unbalanced matches, where games are either too easy or too hard.
Such mismatches are a leading cause of user churn, directly impacting platform revenue and growth.
Thus, the matching policy is not merely a logistical tool.
It is a strategic lever that directly impacts unit economics and long-term viability.

However, managing this trade-off is notoriously difficult because the future is unknown.
The decision to wait is fundamentally a bet that future arrivals will offer matching quality that outweighs the accumulated cost of delay.
In stable environments with predictable arrival patterns, this might be a calculated risk that can be optimized through standard stochastic control methods.
But real-world platforms operate under highly non-stationary conditions: demand spikes unpredictably during rainstorms, driver supply fluctuates with traffic conditions, and player logins fluctuate significantly with time of day and promotional events.
In such volatile environments, the optimal waiting window becomes a moving target that shifts continuously with changing market conditions.
A static rule that works well on average, such as waiting a fixed two minutes, can be disastrous during a demand surge, leading to exploding queue lengths and system overload. 
Conversely, the same rule can be inefficient during a lull, causing unnecessary delays when immediate matching would be optimal.

Despite the practical significance of this problem, the existing toolkit for managing it remains divided between two extremes.
On one side, industry practitioners often rely on ad hoc heuristics, such as fixed time windows (e.g., ``batch every 30 seconds'') or simple queue-length thresholds.
While operationally simple, these static rules lack theoretical robustness and struggle to adapt to the rapid, non-stationary shifts in demand intensity that characterize real-time platforms.
When market conditions change, these heuristics often require manual retuning, creating operational overhead and introducing the risk of suboptimal performance during transition periods.
On the other side, the academic literature offers sophisticated dynamic programming and stochastic control models.
These approaches, while theoretically rigorous, typically depend on strong distributional assumptions (such as known Poisson arrival rates) or require complex state-dependent policies that are computationally brittle in volatile, real-world environments.
Consequently, there is a need for timing rules that retain the operational simplicity of heuristics while offering explicit robustness guarantees in settings where arrival forecasts are unreliable.

We address this need through a queue-length model and a novel operational principle: Cost-Balancing (CB).
Rather than optimizing for a specific, predicted future, we seek an operational equilibrium where the system self-regulates based on realized costs.
Our core insight is that a matching decision should be triggered when the accumulated cost of waiting reaches a calibrated proportion of the instantaneous cost of matching.
This mechanism acts as a state-dependent control lever that adapts to changing market conditions. 
In low-demand periods, the system naturally waits longer to aggregate quality and exploit economies of scale. 
In high-demand spikes, the waiting costs accumulate faster, triggering quicker matches to prevent backlog and system overload.
The algorithm's decision boundary shifts dynamically with the state of the system, so the effective matching threshold responds to whether the market is thick or thin.
This intuition mirrors fundamental economic principles of marginal analysis but is applied here as a robust timing rule that uses realized waiting and matching costs rather than a fitted arrival-distribution model.

\subsection*{Our Contributions}

Our work makes the following contributions to the theory and practice of dynamic matching:

\begin{enumerate}
    \item 
    We formulate a multi-sided matching model that captures general, state-dependent cost structures and non-stationary arrival dynamics, accommodating applications such as ride-pooling, team formation, and multi-resource allocation. Building on this framework, we introduce the Cost-Balancing algorithm, a simple, adaptive policy that implements our equilibrium principle without requiring distributional information or extensive parameter tuning, making it immediately deployable in practice.

    \item 
    We characterize the structural properties of the optimal policy.
    We prove that while the general problem is intractable, a relaxed model with Poisson arrivals and a convex, supermodular matching cost function admits a monotone optimal matching region.
    Using a steady-state fluid approximation of a symmetric benchmark, we derive the cost-balance condition that motivates the ratio form of CB.

    \item 
    We prove that CB is 2-competitive on finite balanced adversarial instances under the robust calibration $\alpha=1$, and more generally that $CB_\alpha$ has competitive ratio $1+\max\{\alpha,1/\alpha\}$.
    This robust guarantee contrasts sharply with standard heuristics: greedy and fixed-threshold policies can incur unbounded costs in worst-case scenarios.
    We further establish a matching lower bound of $2$ over the full class of monotone complete-tuple matching-cost functions, showing that $CB_1$ is worst-case optimal among deterministic online algorithms in the base queue-length model.
    We also show that this factor-$2$ guarantee extends to any single fixed heterogeneous consumption requirement.

    \item 
    We demonstrate the algorithm's practical effectiveness through extensive experiments on game matchmaking and real-world food delivery data.
    A key advantage of our approach is its state-adjusting trigger: unlike fixed time or quantity thresholds, the CB mechanism changes its effective threshold with realized waiting costs, matching costs, and pool sizes.
    Our algorithm outperforms industry-standard batching heuristics across diverse market conditions.
\end{enumerate}

The remainder of the paper is organized as follows.
Section~\ref{sec:literature} reviews the related literature, positioning our work within three research streams.
Section~\ref{sec:model} formally defines the multi-sided matching model, introduces the Cost-Balancing algorithm, and states our main result on competitive ratio guarantees.
Section~\ref{sec:principle} develops the theoretical foundation by analyzing the structure of optimal policies in a relaxed model and deriving a cost-balance benchmark from a steady-state fluid approximation.
Section~\ref{sec:analysis} presents the competitive ratio analysis of the Cost-Balancing algorithm, demonstrates the fragility of standard benchmark policies, and establishes the optimality of the factor-$2$ guarantee.
Section~\ref{sec:extensions} develops the fixed heterogeneous consumption extension.
Section~\ref{sec:numerical} evaluates the algorithm through calibrated numerical experiments on game matchmaking and food delivery data, providing evidence of practical effectiveness across diverse market conditions.
Finally, Section~\ref{sec:conclusion} concludes the paper and suggests directions for future research.

\section{Literature Review}
\label{sec:literature}

Dynamic matching has emerged as an important research domain bridging operations research, algorithmic economics, and theoretical computer science.
This paper proposes a cost-balancing mechanism that formalizes the temporal equilibrium between immediate matching costs and the opportunity cost of delay.
Our contribution engages with three literature streams: (1) competitive ratio analysis of online matching algorithms, (2) the matching-waiting trade-off in operations, and (3) the effectiveness of simple mechanisms in complex environments.

\subsection{Competitive Ratio Analysis}

The theoretical computer science community has pioneered competitive ratio analysis for online bipartite matching.
The seminal work of \cite{KVV1990} established a competitive ratio of $1-\frac{1}{e}$ for adversarial arrivals, spawning two key extensions: stochastic arrival models \citep{FMMM2009, BK2010, MGS2012, JL2014} that derive improved ratios under distributional assumptions, and stochastic reward frameworks \citep{MP2012, HZ2020, GU2023} that extend the analysis to uncertain match values; see \cite{Mehta2013} for a comprehensive survey.
Related beyond-worst-case work on online metric matching uses stochastic arrivals or predictions to improve performance guarantees \citep{YY2026}.
Our work contributes a constant-competitive algorithm to this literature that requires no distributional assumptions and extends from bipartite to multi-sided matching.

\subsection{The Matching-Waiting Trade-off}

A central challenge in dynamic matching systems, extensively studied within the operations research community, lies in balancing the immediate rewards of matching against the long-term optimization potential afforded by waiting.
Our modeling framework is closely related to the matching queues literature \citep{GW2015, ACG2022}.
\cite{GW2015} analyze the dynamic control of matching queues to minimize holding costs, and \cite{ACG2022} explore the optimal design of matching topologies to balance rewards and delays.

The contemporary literature formalizes this trade-off through distinct optimization frameworks.
Multi-stage optimization approaches, such as those in \cite{FNS2024} and \cite{FN2025}, achieve competitive vertex-weighted matching by dynamically transitioning between greedy and hedging strategies.
A complementary line of work \citep{WXY2023, KAG2024, KAG2025, CKKZ2026} adopts regret minimization to benchmark policies against optimal hindsight decisions, often in settings with heterogeneous demand and supply.

Platform-specific dynamics further shape this trade-off, particularly in spatial settings where physical constraints impose unique challenges.
\cite{kanoria2025} characterize the scaling behavior of achievable costs in dynamic spatial matching with uncertain locations, while related spatial models study staffing and service-range design in on-demand platforms \citep{BH2026, ASY2026}.
In ride-sharing, \cite{ABDJSS2019} and \cite{AS2020} incorporate stochastic abandonment and heterogeneous sojourn times, while \cite{WZZ2024} model abandonment and cancellation dynamics, and \cite{LWX2025} derive exact cost formulas for spatial matching on a circle.
Recent work on impatient and heterogeneous demand and supply similarly emphasizes the need to balance match value against abandonment risk \citep{ADSW2025}.
In on-demand delivery, \cite{GG2024} propose a cost-based threshold policy for delivery dispatch, while \cite{CH2024} and \cite{MMR2025} focus on delay-sensitive dispatching and spatial pooling, respectively.
More broadly, the link between market thickness and matching quality has been explored across various platform contexts \citep{AR2021, HZ2022, ILMW2023, KSS2024, ZPT2024, CEL2026}, with \cite{KG2025} specifically quantifying the ``cost of impatience'' via scaling laws.

A related stream of literature investigates the mechanics of batching and delayed decisions.
\cite{XMX2026} and \cite{WKJ2021} study the benefits of batching in online matching, and \cite{BEV2026} examine optimal batching schedules with competitive guarantees.
Our work differs from these contributions by studying a realized-cost timing rule for multi-type matching systems and by providing a distribution-free competitive guarantee in the complete-tuple and fixed-consumption settings.

\subsection{The Effectiveness of Simple Policies}

A growing body of literature in economics and operations management demonstrates that simple heuristic policies can achieve near-optimal performance in dynamic matching markets \citep{ALG2020, BLY2020, MNP2020, LMT2022, ANS2023}.
Particularly relevant is \cite{MNP2020}, who study the dynamic clearing game in two-sided markets with Poisson arrivals and characterize the conditions under which simple policies are optimal.

Subsequent research explores various simple heuristics across different settings.
\cite{BRSWW2022} analyze departure-based threshold policies under general random utilities, while \cite{BFP2023} and \cite{Gupta2024} establish performance bounds for greedy algorithms on specific graph structures.
\cite{ET2026} show that simple batching policies achieve near-optimal performance even with impatient agents.
\cite{KAG2024} introduce the ``general position gap'' to characterize when periodic clearing approximates optimal values in multiway matching, and their follow-up work \citep{KAG2025} proves that greedy heuristics can be hindsight optimal in two-way settings.
\cite{MMR2025} propose a potential-based heuristic for delivery pooling under specific reward structures.

A common thread in this literature is the reliance on stochastic assumptions to establish performance guarantees.
Our Cost-Balancing algorithm achieves a 2-competitive guarantee for finite balanced adversarial arrival sequences under the robust calibration $\alpha=1$, without distributional requirements or asymptotic market conditions.
Despite its simple design, our policy has a stronger worst-case guarantee than standard benchmark heuristics and performs well in the calibrated experiments, highlighting its robustness and practical value.

\section{Model, Algorithm, and Main Results}
\label{sec:model}

In this section, we develop a general framework for multi-sided matching that captures the essential structural elements of these markets: stochastic arrivals, the cost of waiting, and the economies of scale inherent in matching. 
This section formalizes the problem environment, introduces the Cost-Balancing algorithm as a robust solution mechanism, and states our main performance guarantee.

\subsection{The Multi-sided Matching Problem}
\label{sec:problem definition}

We consider a continuous-time matching system with $N \geq 2$ distinct types of agents.
Agents of each type $i \in \{1, \ldots, N\}$ arrive dynamically over time and must be matched into complete $N$-tuples (one agent from each type) to complete service.
Let $A_i(t)$ denote the counting process of type-$i$ arrivals by time $t$, with arrival epochs $\{\sigma_{i,k}\}_{k \geq 1}$ defined by $\sigma_{i,0} = 0$ and $\sigma_{i,k} = \inf\{t : A_i(t) \geq k\}$ for $k \geq 1$.
We impose no distributional assumptions on these arrival streams: they may be deterministic, stochastic, or non-stationary, and are assumed only to be non-explosive (i.e., finitely many arrivals in any finite interval).
This generality allows our framework to accommodate the volatile and potentially adversarial arrival patterns characteristic of real-world platforms.

\begin{figure}[htbp]
    \centering
    \begin{tikzpicture}[
        node distance=0.8cm,
        box/.style={draw, minimum width=2.2cm, minimum height=0.6cm, 
        font=\small, align=center},
        decision/.style={diamond, draw, fill=gray!10, aspect=1.5}]
        
    \node[box] (arrive) {Arrival Streams};
    \node[box, below=of arrive] (queue) {Type Queues \\ $X_1,\dots,
    X_N$};
    \node[decision, below=of queue] (decide) {Match Algorithm};
    \node[box, below left=1cm of decide] (wait) {Back to Queues \\ 
    (Waiting Cost $\uparrow$)};
    \node[box, below right=1cm of decide] (exit) {Tuple Leave \\ 
    (Matching Cost)};
    
    \draw[->, thick] (arrive) -- (queue);
    \draw[->, thick] (queue) -- (decide);
    \draw[->, thick] (decide) -- node[left] {\footnotesize Wait} (wait);
    \draw[->, thick] (decide) -- node[right] {\footnotesize Match} 
    (exit);
    \draw[->, thick, dashed] (wait) -- ++(-2,0) |- (queue);
    \end{tikzpicture}
    \caption{System Dynamics and Decision Flow}
    \label{fig:flowchart}
\end{figure}
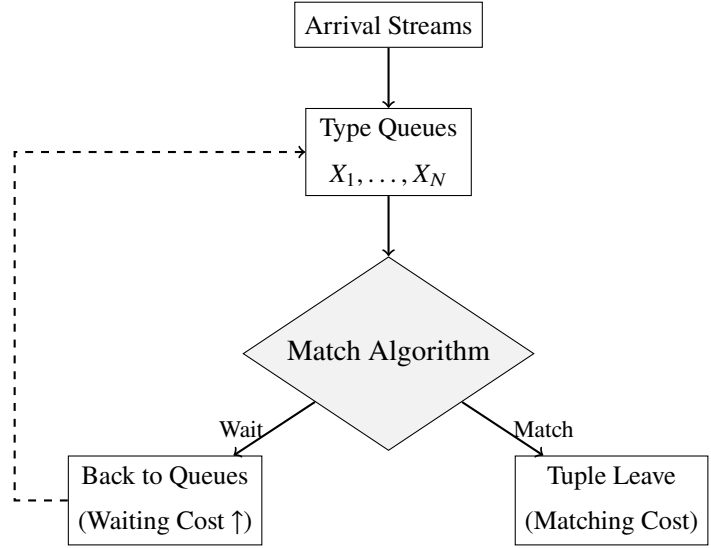

The state of the system at time $t \geq 0$ is captured by the vector of queue lengths $\mathbf{X}(t) = (X_1(t), \ldots, X_N(t)) \in \mathbb{Z}_+^N$, where $X_i(t)$ denotes the number of unmatched type-$i$ agents.
Let
\[
\mathcal F := \{\mathbf{x}\in\mathbb Z_+^N:\min_i x_i\geq 1\}
\]
denote the feasible region.
A matching decision is feasible at time $t$ if and only if $\mathbf{X}(t)\in\mathcal F$, i.e., at least one agent of each type is available.
The controller may execute a match at any calendar time.
The state changes only at arrivals and matches, so an implementable policy can be represented through arrival events, matching events, and policy-generated timer events.

We use the following event-ordering convention throughout the competitive analysis.
Events are processed one at a time.
If multiple arrivals occur at the same calendar time, they are processed in an arbitrary but fixed order, and the policy may re-evaluate its matching decision after each individual arrival.
If a timer event and an arrival occur at the same calendar time under the CB policy, the timer event is processed first.
Upon arrival of a type-$i$ agent, the queue length increases by one, $\mathbf{X}\leftarrow \mathbf{X}+\mathbf e_i$, where $\mathbf{e}_i$ is the $i$-th standard basis vector.
If a match is executed, one agent from each type is selected to form an $N$-tuple that exits the system, and the state transitions as $\mathbf{X}\leftarrow \mathbf{X}-\mathbf 1$, where $\mathbf{1}=(1,\ldots,1)^\top\in\mathbb Z_+^N$.

For any policy $\pi$, let $\tau_k^\pi$ denote the calendar time of its $k$-th match.
Let $\mathbf{Y}_k^\pi$ denote the state on which the $k$-th matching decision is made: after all arrivals assigned to that event have been processed and before the matched tuple leaves the system.
Each executed match incurs a state-dependent cost $f(\mathbf{Y}_k^\pi)$, where $f: \mathbb{Z}_+^N \to \mathbb{R}_+$ represents the minimum cost of forming a matched tuple given the current pool of available agents.
This cost captures factors such as geographic dispersion in ride-sharing, skill disparity in gaming, or delivery distance in logistics.
We impose the following structural assumption on the matching cost function. Note that $\mathbf{x}\succeq \mathbf{y}$ means $\mathbf{x}[i]\ge \mathbf{y}[i]$ for each $i$, where $\mathbf{x}[i]$ and $\mathbf{y}[i]$ denote the $i$-th components of $\mathbf{x}$ and $\mathbf{y}$, respectively.

\begin{assumption}
    \label{ass:key}
    The matching cost incurred by a feasible match is a function $f(\cdot)$ of the pre-match queue-length vector.
    The function is positive and finite on $\mathcal F$: $0<f(\mathbf{x})<\infty$ for all $\mathbf{x}\in\mathcal F$.
    Furthermore, $f(\cdot)$ satisfies the monotonicity property on feasible states, i.e., $f(\mathbf{x})\le f(\mathbf{y})$ whenever $\mathbf{x},\mathbf{y}\in\mathcal F$ and $\mathbf{x} \succeq \mathbf{y}$.
\end{assumption}

Assumption~\ref{ass:key} is the maintained abstraction for the main theoretical analysis; it captures the market thickness effect: as the pool of available agents grows, the platform can select better matches, thereby reducing matching costs.
This property arises naturally in many applications.
In ride-sharing, larger driver pools enable better geographic clustering; in gaming, more players allow tighter skill matching; in delivery, more orders permit efficient route bundling.
We adopt a queue-dependent cost structure rather than modeling individual agent attributes explicitly.
This modeling choice is both practically motivated and theoretically necessary.
In practice, queue length serves as a sufficient statistic that aggregates the benefit of agent heterogeneity.
Theoretically, Appendix~\ref{sec:general_model} documents a boundary of this abstraction: when the model keeps full attribute information and the adversary can choose attributes, computing optimal policies is NP-hard even in static settings, and no deterministic online algorithm can achieve a bounded competitive ratio under adversarial attribute arrivals.
These results motivate Assumption~\ref{ass:key} as a tractable abstraction that enables meaningful performance guarantees.

Unmatched agents incur waiting costs while residing in their respective queues.
We assume a linear waiting cost structure: each unmatched type-$i$ agent incurs a cost at rate $c_i > 0$ per unit time.
The instantaneous waiting cost rate of the system is defined as
$$
w(\mathbf{X}(t)) := \sum_{i=1}^N c_i X_i(t),
$$
representing the aggregate cost rate across all waiting agents.
The waiting cost parameters $\{c_i\}$ may differ across types, reflecting heterogeneous time-sensitivity (e.g., passengers may be more impatient than drivers).

A policy $\pi$ specifies when to execute matches based on the arrival history, current state, and any policy-generated timer events.
The total cost under policy $\pi$ over a horizon $[0, T]$ consists of cumulative waiting costs and matching costs:
$$
J_T^{\pi} := \int_0^T w(\mathbf{X}^{\pi}(t)) \, \mathrm{d}t + \sum_{k=1}^{R^\pi(T)} f(\mathbf{Y}_k^\pi),
$$
where $R^\pi(T)$ denotes the number of matches executed by policy $\pi$ by time $T$.
For the finite-instance competitive analysis, we focus on balanced complete-tuple instances.
A finite arrival sequence is balanced if there exists an integer $K$ such that exactly $K$ agents of each type arrive.
Since each base-model match consumes one agent of every type, such an instance contains exactly $K$ feasible complete tuples.
Every admissible policy must eventually execute these $K$ complete-tuple matches; if a policy fails to clear all $K$ tuples, its cost is defined to be infinite.
Otherwise, its finite-instance total cost is
\[
J^{\pi}
:=
\int_0^{\tau_K^\pi} w(\mathbf{X}^{\pi}(t))\,\mathrm{d}t
+
\sum_{k=1}^K f(\mathbf{Y}_k^\pi).
\]
The balanced-instance requirement is a scope condition for the finite competitive analysis rather than a claim that real arrival streams are exactly balanced.
On a finite sample with unequal type counts, agents from overrepresented types cannot be cleared by any complete-tuple timing policy unless the model adds another operational convention, such as carryover to a future horizon, abandonment, outsourcing, dummy completion, partial matching, or a terminal penalty.
These conventions introduce additional decisions and cost terms that are separate from the timing trade-off studied in the base model.
We therefore use balanced finite instances to compare policies on a common feasible region: both the online policy and the offline benchmark must clear the same $K$ complete tuples.
For infinite-horizon problems, we consider the long-run average cost:
$$
\bar J^{\pi} := \limsup_{T \to \infty} \frac{1}{T} J_T^{\pi}.
$$
The platform's objective is to find a policy $\pi^*$ that minimizes this cost.
Theorem~\ref{thm:main} below is a finite-instance competitive guarantee; the long-run objective provides the broader modeling benchmark for stochastic settings.

\subsection{The Cost-Balancing Algorithm}
\label{sec:algorithm}

Standard approaches to this problem typically fall into two categories: static policies that rely on fixed parameters (e.g., matching every 5 minutes or when the queue reaches a threshold of 10), and dynamic approaches that require precise knowledge of future arrival rates. 
Both approaches suffer from significant drawbacks in volatile environments. 
Static policies fail to adapt to demand surges, while dynamic policies are sensitive to estimation errors in arrival rate forecasts. 
To overcome these limitations, we introduce the Cost-Balancing algorithm, presented in Algorithm~\ref{alg:CB}.
In the algorithm, a future arrival remains if at least one arrival epoch in the input stream has not yet occurred.

\begin{algorithm}[htbp]
\caption{Cost-Balancing Algorithm}
\label{alg:CB}
\begin{algorithmic}[1]
\Require Arrival streams for each type $i\in\{1,\ldots,N\}$; balancing parameter $\alpha>0$.
\Ensure Total cost $J$.
\State Initialize $\mathbf{x}\gets\mathbf 0$, $W\gets0$, and $J\gets0$.
\While{a future arrival remains or $\mathbf{x}\in\mathcal F$}
    \If{$\mathbf{x}\notin\mathcal F$}
        \State Wait until the next arrival.
        \State Update $W$ over the elapsed time.
        \State Process that arriving type-$i$ agent: $\mathbf{x}\gets\mathbf{x}+\mathbf e_i$.
    \ElsIf{$f(\mathbf{x})\leq\alpha W$}
        \State $J\gets J+W+f(\mathbf{x})$.
        \State Match one complete tuple: $\mathbf{x}\gets\mathbf{x}-\mathbf 1$.
        \State $W\gets0$.
    \Else
        \State Set $\Delta(\mathbf{x},W)=\frac{f(\mathbf{x})/\alpha-W}{w(\mathbf{x})}$.
        \State Wait until the earlier of the next arrival, if any, and the timer expiration after $\Delta(\mathbf{x},W)$.
        \State Update $W$ over the elapsed time.
        \If{the timer expires first}
            \State $J\gets J+W+f(\mathbf{x})$.
            \State Match one complete tuple: $\mathbf{x}\gets\mathbf{x}-\mathbf 1$.
            \State $W\gets0$.
        \Else
            \State Process the arriving type-$i$ agent: $\mathbf{x}\gets\mathbf{x}+\mathbf e_i$.
        \EndIf
    \EndIf
\EndWhile
\State \Return $J$.
\end{algorithmic}
\end{algorithm}

The CB algorithm operates on a simple but powerful principle: trigger a match when the accumulated waiting cost becomes comparable to the instantaneous matching cost.
CB keeps only the queue vector $\mathbf{x}$ and the accumulated waiting cost $W$ since the previous match.
The rule is evaluated whenever an arrival occurs or a timer expires.
If the pool is feasible and
\[
    f(\mathbf{x})\leq \alpha W,
\]
CB clears one complete tuple and resets $W$ to zero, where $\alpha > 0$ is the balancing parameter defining the CB policy family.
If the pool is feasible but the trigger is not yet met, CB sets the timer to the exact time at which the trigger would bind absent further arrivals.
If an arrival occurs before the timer, the timer is discarded and the rule is evaluated again.
Thus, a match can be triggered either by an arrival that lowers the matching cost enough or by the endogenous timer as waiting cost accumulates.

The intuition behind this design is as follows.
Under Assumption~\ref{ass:key}, the matching cost $f(\mathbf{x})$ decreases as the queue lengths grow, reflecting the benefit of market thickness.
Meanwhile, the waiting cost $W$ accumulates over time as agents wait in the system.
These two forces create opposing pressures: waiting longer reduces matching costs but increases waiting costs.
The cost-balancing condition $f(\mathbf{x}) \leq \alpha W$ provides a realized-cost stopping rule for when accumulated delay has become large relative to the current matching cost.
It should be interpreted as a robust proxy rather than as an exact marginal optimality condition.

The CB algorithm offers three key advantages over existing approaches.
First, it is adaptive: unlike fixed thresholds that must be recalibrated when demand patterns change, the CB algorithm adjusts its matching frequency based on realized costs once $\alpha$ is fixed or calibrated.
In high-demand periods, waiting costs accumulate rapidly, triggering more frequent matches; in low-demand periods, the algorithm naturally waits longer to exploit market thickness.
Second, it is computationally simple: implementation requires only tracking the queue lengths and a running sum of waiting costs, with no need for demand forecasting or complex optimization.
Third, as we establish in the following subsection, the algorithm provides provable worst-case guarantees for finite balanced adversarial arrival sequences.

\subsection{Main Result}
\label{sec:main_result}

To evaluate the performance of online matching algorithms, we adopt the competitive ratio framework, which is the standard metric in online decision-making for comparing an algorithm's performance against an omniscient offline benchmark.
Let $J^{OPT}$ denote the minimum total cost achievable by an offline algorithm with complete knowledge of the arrival sequence in advance.
The competitive ratio measures the worst-case performance gap between an online algorithm and this optimal benchmark.

\begin{definition}[Competitive Ratio]
    An online policy $\pi$ is $\rho$-competitive if, for every finite instance and arrival sequence in the model setting,
    \[
        J^\pi \leq \rho J^{OPT}.
    \]
    The competitive ratio of $\pi$ is the infimum over all $\rho$ for which this inequality holds.
\end{definition}

In the base complete-tuple analysis, the relevant finite instances are the balanced arrival sequences satisfying Assumption~\ref{ass:key}, with the finite-instance cost defined above.
A lower competitive ratio indicates better worst-case performance.
An algorithm with competitive ratio $\rho$ guarantees that its cost is at most $\rho$ times the optimal offline cost across the relevant adversarial arrival sequences.
This worst-case framework is particularly relevant for matching platforms, which face demand shocks, seasonal fluctuations, and non-stationary arrival patterns that are difficult to forecast.
A bounded competitive ratio ensures that the algorithm performs well even in scenarios that were not anticipated during design.

Our main theoretical result establishes that the CB algorithm achieves a bounded competitive ratio on finite balanced instances.

\begin{theorem}
    \label{thm:main}
    Consider any finite balanced complete-tuple instance with $K$ arrivals of each type.
    Suppose Assumption~\ref{ass:key} holds.
    For any $\alpha>0$, the Cost-Balancing algorithm $CB_\alpha$ satisfies
    \[
    J^{CB_\alpha}
    \leq
    \left(1+\max\left\{\alpha,\frac{1}{\alpha}\right\}\right)J^{OPT}.
    \]
    In particular, choosing $\alpha=1$ gives
    \[
    J^{CB_1}\leq 2J^{OPT}.
    \]
\end{theorem}

The theorem gives the full tradeoff induced by the balancing parameter.
The term $1+\alpha$ controls cases in which OPT's corresponding match occurs no earlier than CB's match, while the term $1+1/\alpha$ controls cases in which OPT matches earlier.
The robust worst-case calibration is therefore $\alpha=1$, which yields a 2-competitive policy.

Theorem~\ref{thm:main} has several important implications.
First, the CB algorithm provides a guaranteed performance bound: its cost is at most twice the offline optimum when $\alpha=1$.
This guarantee holds for adversarial finite balanced arrival sequences.
Second, the competitive ratio is distribution-free and depends only on the chosen balancing parameter, not on the arrival timings.
Thus, the theorem does not require a stochastic arrival model.
Third, as we will demonstrate in Section~\ref{sec:analysis}, standard heuristics such as greedy matching and fixed-threshold policies can have unbounded competitive ratios, meaning their costs can be arbitrarily worse than the optimum in adversarial scenarios.
The CB algorithm thus offers a qualitatively stronger guarantee than these common approaches.

The proof of Theorem~\ref{thm:main} and additional theoretical results, including a lower bound on the competitive ratio achievable by any online algorithm, are presented in Section~\ref{sec:analysis}.

\section{The Cost-Balancing Principle}
\label{sec:principle}

The Cost-Balancing algorithm introduced in Section~\ref{sec:algorithm} is not an ad hoc heuristic; its design is motivated by structural properties of optimal policies and by a deterministic approximation.
This section develops the design logic underlying the cost-balancing principle.
Our analysis proceeds in two stages.
First, we examine a relaxed model with stationary arrivals to establish that optimal matching regions are monotone, thereby motivating the search for state-dependent matching rules.
Second, we use a steady-state fluid approximation to derive a transparent relationship between waiting and matching cost rates.
This calculation shows how a cost-minimizing stationary operating point can be characterized by a specific cost ratio.
The Cost-Balancing algorithm uses this ratio logic as a timing rule, while the worst-case guarantee in Theorem~\ref{thm:main} supplies a separate robust calibration.

\subsection{Monotonicity of Optimal Policies in Relaxed Settings}
\label{subsec:relaxed_optimality}

A fundamental question in matching algorithm design is whether simple monotone matching rules can be optimal.
In general settings with arbitrary arrival patterns, the answer is often negative: the optimal decision rule can be highly complex and non-monotonic.
However, when we restrict attention to stationary environments with well-behaved cost structures, a clear structure emerges.

We consider a relaxed model with Poisson arrivals and impose additional conditions on the matching cost function.

\begin{definition}
    A matching cost function $f: \mathbb{Z}_+^N \to \mathbb{R}_+$ is called convex, supermodular, and component-wise non-increasing if it satisfies the following properties:
    \begin{itemize}
        \item Convexity: $2f(\mathbf{x} + \mathbf{e}_i) \leq f(\mathbf{x}) + f(\mathbf{x} + 2\mathbf{e}_i)$ for all $\mathbf{x}$ and $i$.
        \item Supermodularity: $f(\mathbf{x} + \mathbf{e}_i) + f(\mathbf{x} + \mathbf{e}_j) \leq f(\mathbf{x}) + f(\mathbf{x} + \mathbf{e}_i + \mathbf{e}_j)$ for all $\mathbf{x}$ and $i \neq j$.
        \item Component-wise non-increasing: $f(\mathbf{x} + \mathbf{e}_i) \leq f(\mathbf{x})$ for all $\mathbf{x}$ and $i$.
    \end{itemize}
\end{definition}

Convexity ensures diminishing marginal returns from additional agents, supermodularity captures complementarity across agent types, and the non-increasing property reflects economies of scale.
Many natural matching cost functions satisfy these conditions.
For example, in spatial matching where agents are uniformly distributed, the expected minimum distance decreases convexly in the number of agents.

Under these conditions, we establish the following structural result.

\begin{proposition}
\label{prop:relaxed_monotonicity}
    Consider the finite-horizon uniformized relaxation of the Poisson-arrival matching system, where at each decision epoch the controller either waits or, if $\mathbf{x}\in\mathcal F$, executes one complete-tuple match.
    Suppose the matching cost function is convex, supermodular, and component-wise non-increasing.
    Then, at every finite horizon, the set of feasible states at which matching is the optimal action is closed upward in queue lengths: if matching is optimal at state $\mathbf{x}\in\mathcal F$, then matching is also optimal at any state $\mathbf{y}\in\mathcal F$ with $\mathbf{y} \succeq \mathbf{x}$.
\end{proposition}

\begin{remark}
    The proof of Proposition~\ref{prop:relaxed_monotonicity}, detailed in Appendix~\ref{app:proof_relaxed_monotonicity}, requires a novel analytical approach.
    Standard dynamic programming arguments establish monotonicity by propagating convexity or supermodularity of the value function across state transitions.
    However, these properties fail to hold in multi-sided matching due to the coupled state transitions: a single match simultaneously reduces all queue lengths by one ($\mathbf{x} \to \mathbf{x} - \mathbf{1}$), breaking the standard induction structure.
    We overcome this challenge by introducing a diagonal difference function $G_k(\mathbf{x}) = J_k(\mathbf{x}) - J_k(\mathbf{x} - \mathbf{1})$ and establishing specialized bounds on its increments relative to the matching cost.
    This technique enables direct induction on the net matching benefit without requiring global structural properties of the value function.
\end{remark}

Proposition~\ref{prop:relaxed_monotonicity} confirms that in stable environments, the intuition ``more agents make matching more attractive'' holds rigorously.
The matching region is monotone: once matching is optimal at a feasible state, it remains optimal after additional agents arrive in any queues.

However, this result also reveals a fundamental limitation: the location of the optimal matching region depends on the arrival rate.
A static queue-length rule that is optimal for one arrival rate will be suboptimal when the rate changes.
In practice, platforms face volatile demand patterns where arrival rates fluctuate significantly over time (e.g., peak vs. off-peak hours, weather events, special promotions).
A fixed threshold policy, no matter how carefully calibrated, will perform poorly when conditions deviate from the assumed baseline.
This observation motivates the search for an adaptive policy that adjusts its matching criterion in real-time without requiring knowledge of the arrival rate.

\subsection{The Cost-Balancing Principle via a Steady-State Fluid Approximation}
\label{subsec:fluid_insight}

To build intuition for how a matching rule should be calibrated in stable environments, we analyze a tractable deterministic steady-state fluid approximation.
The calculation is intended as a design benchmark rather than a stochastic optimality theorem.
Its key insight is that a cost-minimizing operating point is characterized by a cost ratio condition, which forms the basis for the Cost-Balancing algorithm.

For tractability, we consider a symmetric two-sided benchmark ($N=2$) with mean arrival rate $\frac{\lambda}{2}$ for each type and equal waiting cost rates $c_1 = c_2 = c$.
We focus on the two-sided case for clarity of exposition; we expect a similar cost-ratio characterization to hold for general $N$, though the optimal ratio may depend on the specific form of the cost function in the multi-sided case.
We assume the matching cost follows a power-law form that explicitly captures economies of scale:
\begin{assumption}
    \label{ass:illustrative}
    The matching cost function is $f(x_1, x_2) = \frac
    {\kappa}{(x_1 x_2)^{\beta}}$, where $\kappa > 0$ is a scale parameter and $\beta > 0$ controls the strength of economies of scale.
\end{assumption}
This functional form satisfies Assumption~\ref{ass:key}: matching costs decrease as queue lengths grow.
Such power-law cost structures arise naturally in spatial matching settings and have been widely adopted in the literature (e.g., \cite{WZZ2024}).
A larger $\beta$ corresponds to stronger economies of scale, where the benefit of waiting (i.e., reduced matching cost) is more pronounced.

The following proposition states the resulting cost ratio in this deterministic benchmark.

\begin{proposition}
    \label{prop:fluid_balance}
    Consider the symmetric two-sided steady-state fluid benchmark under Assumption~\ref{ass:illustrative}, with mean arrival rate $\lambda/2$ for each type and equal waiting cost rates $c_1=c_2=c$.
    For a common positive fluid queue level $\bar X>0$, the associated steady-state cost rate is
    $$
    \bar J(\bar X)
    =
    2c\bar X
    +
    \frac{\lambda}{2}\frac{\kappa}{\bar X^{2\beta}}.
    $$
    Let $\bar X^*$ be the unique minimizer of $\bar J(\bar X)$ over $\bar X>0$.
    If
    $$
    W_{\mathrm{rate}}^*=2c\bar X^*,
    \qquad
    M_{\mathrm{rate}}^*=\frac{\lambda}{2}\frac{\kappa}{(\bar X^*)^{2\beta}},
    $$
    denote the waiting-cost and matching-cost rates at this benchmark optimum, then
    $$ \frac{W_{\mathrm{rate}}^*}{M_{\mathrm{rate}}^*}=2\beta. $$
\end{proposition}

Proposition~\ref{prop:fluid_balance} provides a simple calibration benchmark.
The ratio $2\beta$, which measures the elasticity of the matching cost with respect to the common queue level, determines the steady-state balance point.
When $\beta$ is large (strong economies of scale), the fluid benchmark tolerates higher waiting-cost rates relative to matching-cost rates, because the benefit of accumulating more agents is substantial.
When $\beta$ is small (weak economies of scale), the benchmark matches more promptly, as waiting provides limited benefit.

This calculation motivates the ratio form of the CB algorithm and suggests how $\alpha$ may be calibrated in stable stochastic environments.
The algorithm triggers a match when the condition $M \leq \alpha W$ is satisfied.
Equivalently, $W/M\ge 1/\alpha$.
Here $W_{\mathrm{rate}}^*$ and $M_{\mathrm{rate}}^*$ are cost rates; normalizing both by the same steady-state match rate gives the same ratio on a per-match basis.
Thus the fluid ratio $W_{\mathrm{rate}}^*/M_{\mathrm{rate}}^*=2\beta$ suggests the environment-specific calibration
$$
\alpha_{\mathrm{fluid}}=\frac{1}{2\beta}.
$$
This fluid calibration is distinct from the robust worst-case calibration $\alpha=1$ used in Theorem~\ref{thm:main}.
The fluid calculation is therefore not a literal derivation of Algorithm~\ref{alg:CB}; it provides a steady-state design benchmark for the cost ratio, while the implemented timing rule is justified separately by the competitive-ratio analysis.

The cost-balancing principle also explains the algorithm's robustness to changing conditions.
Consider what happens when the arrival rate suddenly increases.
Higher arrival rates cause waiting costs to accumulate faster (more agents waiting per unit time), while matching costs decrease faster (queues grow more quickly).
Both effects push the ratio $W/M$ upward, causing the cost-balancing condition to be satisfied sooner.
The algorithm thus automatically matches more frequently during high-demand periods.
Conversely, when arrival rates decrease, waiting costs accumulate slowly, and the algorithm naturally waits longer to exploit market thickness.
This automatic adaptation occurs without any knowledge of the arrival rate, because the algorithm responds directly to realized costs rather than forecasted demand.

Proposition~\ref{prop:fluid_balance} and Theorem~\ref{thm:main} play distinct and complementary roles.
Proposition~\ref{prop:fluid_balance} provides design motivation: it shows that in a symmetric two-sided steady-state fluid approximation with power-law costs, the cost-minimizing operating point satisfies $W_{\mathrm{rate}}^*/M_{\mathrm{rate}}^* = 2\beta$, revealing why a ratio-based trigger is natural.
Theorem~\ref{thm:main} provides a performance guarantee: CB is 2-competitive on finite balanced adversarial instances under $\alpha=1$, and the full family $CB_\alpha$ has competitive ratio $1+\max\{\alpha,1/\alpha\}$.
The key point is that the fluid approximation identifies a useful cost-balance benchmark, while the competitive analysis identifies a robust worst-case calibration.

\section{Competitive Ratio Analysis}
\label{sec:analysis}

This section provides rigorous performance guarantees for the Cost-Balancing algorithm and establishes fundamental limits on what any online algorithm can achieve.
Our analysis yields three main results.
First, we prove that CB is 2-competitive on finite balanced instances under $\alpha=1$, with the full family bound $1+\max\{\alpha,1/\alpha\}$ (Section~\ref{subsec:proof_main}).
Second, we demonstrate that standard heuristics, including greedy and fixed-threshold policies, have unbounded competitive ratios, meaning they can perform arbitrarily worse than the optimum in adversarial scenarios (Section~\ref{subsec:benchmarks}).
Third, we prove a matching factor-$2$ lower bound over the full class of monotone complete-tuple matching-cost functions, so the upper and lower bounds together identify the optimal deterministic competitive ratio in the base model (Section~\ref{subsec:lower_bound}).
These results establish the benchmark guarantee for the complete-tuple queue-length framework defined in Section~\ref{sec:model}; Section~\ref{sec:extensions} then examines how far the guarantee extends beyond that base structure.

\subsection{Proof of the Main Theorem}
\label{subsec:proof_main}

We now prove Theorem~\ref{thm:main}.
The proof uses a rank-based charging argument.
For each $k$, we compare CB's $k$-th match with OPT's $k$-th match.
If OPT's $k$-th match is no earlier than CB's, we charge CB's cost to OPT's waiting cost over the same CB intermatch interval.
If OPT's $k$-th match is earlier, we charge CB's cost to OPT's $k$-th matching cost.

\proof{Proof of Theorem~\ref{thm:main}.}
    Fix $\alpha>0$.
    Let $\tau_k=\tau_k^{CB_\alpha}$ be the time of CB's $k$-th match, and set $\tau_0=0$.
    Define the $k$-th CB intermatch interval by $I_k=(\tau_{k-1},\tau_k]$ and let
    \[
        W_k:=\int_{\tau_{k-1}}^{\tau_k} w(\mathbf{X}^{CB_\alpha}(t))\,\mathrm{d}t
    \]
    be CB's waiting cost accumulated during this interval.
    Let $\mathbf{x}_k=\mathbf{Y}_k^{CB_\alpha}$ be CB's pre-match state for its $k$-th match, and define $M_k=f(\mathbf{x}_k)$.
    CB's cost associated with its $k$-th match is
    \[
        C_k^{CB_\alpha}=W_k+M_k.
    \]
    By the CB rule, every CB match occurs only when
    \[
        M_k\leq \alpha W_k.
    \]
    Therefore,
    \begin{equation}
        \label{eq:cb-upper-trigger}
        C_k^{CB_\alpha}\leq (1+\alpha)W_k.
    \end{equation}

    Let $s_k$ be the time of OPT's $k$-th match, and let $\mathbf{u}_k=\mathbf{Y}_k^{OPT}$ be OPT's pre-match state for that match.
    We split the proof into two cases.
    
    First suppose $s_k\geq \tau_k$.
    For almost every $t\in I_k$, CB has completed exactly $k-1$ matches, while OPT has completed at most $k-1$ matches.
    Let $\mathbf A(t)=(A_1(t),\ldots,A_N(t))$ be the cumulative arrival vector.
    Then, outside event times of measure zero,
    \[
        \mathbf{X}^{CB_\alpha}(t)=\mathbf A(t)-(k-1)\mathbf 1,
    \]
    whereas
    \[
        \mathbf{X}^{OPT}(t)=\mathbf A(t)-m^{OPT}(t)\mathbf 1,
        \qquad m^{OPT}(t)\leq k-1,
    \]
    where $m^{OPT}(t)$ is the number of matches completed by OPT by time $t$.
    Hence $\mathbf{X}^{OPT}(t)\succeq \mathbf{X}^{CB_\alpha}(t)$ and, because waiting costs are positive and linear,
    \[
        w(\mathbf{X}^{OPT}(t))\geq w(\mathbf{X}^{CB_\alpha}(t)).
    \]
    It follows that
    \[
        \int_{I_k}w(\mathbf{X}^{OPT}(t))\,\mathrm{d}t\geq W_k.
    \]
    Combining this inequality with \eqref{eq:cb-upper-trigger} yields
    \begin{equation}
        \label{eq:late-charge}
        C_k^{CB_\alpha}
        \leq
        (1+\alpha)\int_{I_k}w(\mathbf{X}^{OPT}(t))\,\mathrm{d}t.
    \end{equation}
    
    Now suppose $s_k<\tau_k$.
    Immediately before their respective $k$-th matches, both policies have completed $k-1$ matches.
    Because $s_k<\tau_k$, the cumulative arrival vector processed before OPT's $k$-th matching decision is componentwise no larger than the cumulative arrival vector processed before CB's $k$-th matching decision.
    Therefore,
    \[
        \mathbf{u}_k\preceq \mathbf{x}_k.
    \]
    By Assumption~\ref{ass:key},
    \[
        M_k=f(\mathbf{x}_k)\leq f(\mathbf{u}_k).
    \]
    
    We next show that $W_k\leq f(\mathbf{u}_k)/\alpha$.
    If CB's $k$-th match is caused by a timer event, then the timer fires exactly when
    \[
        M_k=\alpha W_k.
    \]
    Since $M_k\leq f(\mathbf{u}_k)$, it follows that
    \[
        W_k\leq \frac{1}{\alpha}f(\mathbf{u}_k).
    \]
    
    If CB's $k$-th match is caused by an arrival, let $\mathbf{y}_k$ be CB's state immediately before the triggering arrival is processed.
    We first claim that $\mathbf{y}_k\in\mathcal F$.
    If not, then for some type $i$, the number of type-$i$ arrivals processed immediately before the triggering arrival is at most $k-1$, because CB has removed exactly $k-1$ type-$i$ agents before its $k$-th match.
    Since $s_k<\tau_k$, OPT's $k$-th match was decided with no more type-$i$ arrivals than this, contradicting feasibility of OPT's $k$-th complete-tuple match.
    Hence $\mathbf{y}_k\in\mathcal F$.
    Since CB did not match immediately before the triggering arrival, and since a timer event is processed before an arrival when both occur at the same calendar time, the CB threshold could not have been met just before the arrival:
    \[
        f(\mathbf{y}_k)>\alpha W_k.
    \]
    Moreover, immediately before OPT's $k$-th match and immediately before CB's triggering arrival, both policies have completed $k-1$ matches, and the latter event occurs later.
    Thus $\mathbf{u}_k\preceq \mathbf{y}_k$.
    By Assumption~\ref{ass:key},
    \[
        f(\mathbf{u}_k)\geq f(\mathbf{y}_k)>\alpha W_k,
    \]
    so again $W_k\leq f(\mathbf{u}_k)/\alpha$.
    
    Combining $M_k\leq f(\mathbf{u}_k)$ and $W_k\leq f(\mathbf{u}_k)/\alpha$ gives
    \begin{equation}
        \label{eq:early-charge}
        C_k^{CB_\alpha}
        \leq
        \left(1+\frac{1}{\alpha}\right)f(\mathbf{u}_k).
    \end{equation}
    
    Let
    \[
        \mathcal L=\{k:s_k\geq \tau_k\},
        \qquad
        \mathcal E=\{k:s_k<\tau_k\}.
    \]
    Summing \eqref{eq:late-charge} over $k\in\mathcal L$ and \eqref{eq:early-charge} over $k\in\mathcal E$ yields
    \[
        J^{CB_\alpha}
        \leq
        \rho(\alpha)
        \left[
        \sum_{k\in\mathcal L}\int_{I_k}w(\mathbf{X}^{OPT}(t))\,\mathrm{d}t
        +
        \sum_{k\in\mathcal E}f(\mathbf{u}_k)
        \right],
    \]
    where
    \[
        \rho(\alpha)=1+\max\left\{\alpha,\frac{1}{\alpha}\right\}.
    \]
    The intervals $\{I_k:k\in\mathcal L\}$ are disjoint and lie within OPT's clearing horizon, so the first summation charges disjoint portions of OPT's waiting cost.
    The second summation charges distinct OPT matching costs.
    The bracketed term is therefore at most $J^{OPT}$, proving
    \[
        J^{CB_\alpha}\leq \rho(\alpha)J^{OPT}.
    \]
    Choosing $\alpha=1$ balances the two terms in $\rho(\alpha)$ and gives
    \[
        J^{CB_1}\leq 2J^{OPT}.
    \]
\hfill\Halmos
\endproof

\begin{remark}
    The rank-based proof relies on a common consumption vector: after the same number of matches, CB and OPT must have removed the same number of agents from each type.
    In the complete-tuple model this vector is $\mathbf 1$.
    The same argument extends to any single fixed heterogeneous consumption requirement, as shown in Section~\ref{sec:extensions}.
    It does not directly extend to settings with variable clearing actions or compatibility-network matching, where different policies may consume different type vectors before their $k$-th matches.
\end{remark}

\subsection{The Fragility of Standard Heuristics}
\label{subsec:benchmarks}

To appreciate the value of this guarantee, it is instructive to compare it with the performance of standard policies widely used in practice.
We show that while heuristics like Greedy and Fixed Thresholds may perform well in stable, average-case scenarios, they are fundamentally brittle: they lack the ``safety net" mechanism of Cost-Balancing and can incur unbounded relative costs in adversarial environments.

The Greedy policy matches immediately whenever feasible ($\min_i X_i \geq 1$), prioritizing minimal waiting at the expense of matching quality.
While greedy policies can achieve near-optimal performance under specific conditions such as certain graph structures \citep{Gupta2024}, they generally underperform when matching costs exhibit significant economies of scale.

\begin{proposition}
    \label{prop:greedy-competitive-ratio}
    The competitive ratio of the Greedy policy is unbounded.
\end{proposition}

The Threshold policy waits until the queue length reaches a fixed threshold $\theta$ before matching, deliberately accumulating agents to exploit economies of scale.
This approach is prevalent in periodic clearing models \citep{KAG2024} and can achieve asymptotically optimal performance in thick markets \citep{ET2026}.
However, the optimal threshold depends on arrival rates, making fixed thresholds fragile to demand fluctuations.

\begin{proposition}
    \label{prop:threshold-competitive-ratio}
    The competitive ratio of the Threshold policy is unbounded.
\end{proposition}

The unbounded competitive ratios of Greedy and Threshold policies stem from a common structural flaw: both commit to a fixed decision rule that ignores the realized state of the system.
Greedy always matches immediately, regardless of how much matching cost could be saved by waiting; Threshold always waits to $\theta$, regardless of how much waiting cost has accumulated.
In adversarial scenarios, these rigid commitments can be exploited: an adversary can construct arrival sequences that maximize the gap between the policy's cost and the optimal cost.

The CB algorithm avoids this vulnerability by conditioning its decision on the ratio of realized costs.
The cost-balancing condition $M \leq \alpha W$ creates a decision boundary that shifts with the state of the system: the more waiting cost has accumulated, the higher matching cost the algorithm is willing to accept.
Within the finite balanced complete-tuple model, this state-dependent threshold is the key mechanism underlying the bounded competitive ratio guarantee for CB.

\subsection{Optimality of the Factor 2}
\label{subsec:lower_bound}

Having shown that CB achieves a bounded competitive ratio on finite balanced instances while standard heuristics do not, a natural question arises: how close to optimal can any online algorithm be?
We establish a matching lower bound showing that the factor $2$ in Theorem~\ref{thm:main} cannot be improved by any deterministic online algorithm over the full class of monotone complete-tuple matching-cost functions.
This is a worst-case statement over the model class, not a claim about every fixed matching-cost function.
For example, when matching costs are constant, waiting cannot reduce matching cost and immediate matching can be optimal.
The lower-bound construction instead uses a monotone setup-cost-like function for which batching two complete tuples is almost no more expensive than matching one.

\begin{proposition}
    \label{prop:lower-bound}
    For every $\eta>0$ and every deterministic online policy $\pi$, there exists a matching-cost function $f$ satisfying Assumption~\ref{ass:key} and a finite balanced complete-tuple arrival sequence such that
    \[
        \frac{J^\pi}{J^{OPT}}\geq 2-\eta.
    \]
    Consequently, no deterministic online algorithm can achieve a uniform competitive ratio strictly smaller than $2$ over the class of monotone complete-tuple matching-cost functions.
\end{proposition}

Proposition~\ref{prop:lower-bound} shows that the upper bound in Theorem~\ref{thm:main} is tight.
The lower-bound instance creates a setup-cost-like economy of scale: matching a singleton complete tuple costs $1$, while the incremental cost of matching from a pool with at least two complete tuples can be made arbitrarily small.
The adversary then forces every deterministic online algorithm to choose between clearing singleton tuples before the next arrival and waiting long enough that OPT can avoid the delay.
Together with Theorem~\ref{thm:main}, this establishes $2$ as the optimal deterministic competitive ratio for the base monotone complete-tuple model.

\section{Beyond Complete Tuples: Fixed Heterogeneous Consumption}
\label{sec:extensions}

The main competitive-ratio theorem is proved for a complete-tuple model in which every match consumes one agent from each type.
This section asks how far the guarantee extends when this symmetry is relaxed.
This question is operationally relevant because many platforms clear service bundles that require different numbers of different resources.
For example, delivery dispatch may involve drivers, customer orders, and restaurants or stores; ride-pooling may involve vehicles and riders; and field-service operations may combine technicians, equipment, and jobs.
These examples share a common feature: a clearing action may consume resources in unequal proportions, so the one-agent-per-type abstraction can be too restrictive.

We study a tractable version of this issue in which each clearing action consumes the same vector of resources.
We show that the factor-$2$ guarantee survives when each match consumes such a fixed vector of agents, such as two agents of one type and three of another.
The message is that the rank-based proof depends on common resource accounting: CB and OPT remain comparable after the same number of matches only when every feasible match removes the same type vector.

We consider a fixed consumption vector
\[
    \mathbf r=(r_1,\ldots,r_N)\in\mathbb Z_{++}^N.
\]
A match consumes $r_i$ agents of type $i$ (if a type has zero consumption in a broader specification, it is inactive for this submarket and can be removed from the finite-instance analysis).
The feasible region is
\[
    \mathcal F_{\mathbf r}
    :=
    \{\mathbf x\in\mathbb Z_+^N:x_i\ge r_i \text{ for every } i\}.
\]
A finite instance is $\mathbf r$-balanced if there exists an integer $K$ such that the total number of type-$i$ arrivals is $K r_i$ for every type.
Every admissible policy must eventually execute exactly $K$ matches.

Let $f_{\mathbf r}:\mathcal F_{\mathbf r}\to\mathbb R_+$ be the matching-cost function for this fixed consumption requirement.
We assume $f_{\mathbf r}$ is positive, finite, and componentwise non-increasing on $\mathcal F_{\mathbf r}$.
The fixed-consumption version of CB, denoted $CB_\alpha(\mathbf r)$, is Algorithm~\ref{alg:CB} with two changes: feasibility is evaluated using $\mathcal F_{\mathbf r}$, and a match changes the state as $\mathbf X\leftarrow \mathbf X-\mathbf r$.
The trigger remains
\[
    f_{\mathbf r}(\mathbf X(t))\le \alpha W(t),
\]
where $W(t)$ is the accumulated waiting cost since the previous match.
We use the same event-ordering convention as in Section~\ref{sec:model}: if a timer and an arrival occur at the same calendar time, the timer is processed first.

\begin{proposition}
    \label{prop:fixed-consumption}
    Consider any finite $\mathbf r$-balanced instance.
    Suppose $f_{\mathbf r}$ is positive, finite, and componentwise non-increasing on $\mathcal F_{\mathbf r}$.
    For any $\alpha>0$, the fixed-consumption policy $CB_\alpha(\mathbf r)$ satisfies
    \[
        J^{CB_\alpha(\mathbf r)}
        \le
        \left(1+\max\left\{\alpha,\frac{1}{\alpha}\right\}\right)
        J^{OPT(\mathbf r)}.
    \]
    In particular,
    \[
        J^{CB_1(\mathbf r)}\le 2J^{OPT(\mathbf r)}.
    \]
    Moreover, the factor $2$ is tight over the class of positive, componentwise non-increasing fixed-consumption matching-cost functions.
\end{proposition}

Proposition~\ref{prop:fixed-consumption} shows that the factor-$2$ guarantee and its tightness are not artifacts of one-for-one complete-tuple symmetry.
What matters is that all feasible matches consume the same fixed vector of resources.
In such systems, CB and OPT remain comparable by match rank: after $k-1$ matches, both policies have removed $(k-1)\mathbf r$ whenever they have completed the same number of matches.
Thus the economic mechanism is fixed resource accounting, while the lower-bound mechanism remains uncertainty about future market thickness.
This directly covers service actions that require, for example, two agents of one type and three of another, as long as that requirement is fixed over time.
The same condition also marks the limitation of the extension.
The proposition does not cover flexible dispatch settings in which one driver may serve different numbers of customers across dispatches, or in which the platform chooses variable route bundles from the current pool.
In those environments, different policies can consume different type vectors before their $k$-th clearing action, so the rank-based comparison used in the proof no longer applies directly.
Such flexible one-to-many settings are still operationally important, and the delivery experiment in Section~\ref{sec:delivery} studies one empirically, but they remain outside the formal competitive-ratio guarantee.

\section{Numerical Study}
\label{sec:numerical}

This section evaluates the practical effectiveness of the Cost-Balancing algorithm through calibrated experiments on two real-world applications: video game matchmaking and on-demand food delivery.
The two experiments highlight complementary aspects of the framework.
The gaming study evaluates cost balancing in a continuous-attribute matching environment, where matching quality depends on skill ratings.
The delivery study tests whether the same timing principle remains useful in a richer one-to-many dispatch environment with assignment and routing constraints.
In both settings, the policy uses realized waiting and matching costs to decide when the current pool should be cleared; the balancing parameter is selected on calibration data or by grid search for the corresponding implementation.

\subsection{Video Game Matchmaking}
\label{sec:video-game}

We first examine a 1-vs-1 competitive gaming matchmaking system, where the platform must balance two objectives: minimizing player waiting time and ensuring skill-balanced matches.
This setting introduces complexities absent from our theoretical model: agent heterogeneity is continuous (Elo skill ratings) rather than categorical.

We use a calibrated simulation environment based on the gaming setting of \cite{Gan2023}, rather than a direct replication of their parameterization.
Player arrivals follow a non-stationary Poisson process with a random 24-hour arrival-rate profile.
Each arriving player has a skill rating drawn from $\mathcal{N}(325,75^2)$, truncated below at zero.
The platform's objective is to minimize total cost, defined as
$$
J = \sum_{k=1}^{n} w_k + \gamma \sum_{(i,j) \in \mathcal{M}} |s_i - s_j|,
$$
where $n$ is the number of players, $\mathcal{M}$ is the set of matched player pairs, $w_k$ denotes the waiting time of player $k$, $|s_i - s_j|$ is the skill gap for matched pair $(i,j)$, and $\gamma > 0$ is a weight parameter governing the trade-off between responsiveness and match quality.
We report results for $\gamma=0.1,0.2,\ldots,0.9$, covering the main range from waiting-cost-dominant to matching-cost-sensitive regimes.
The matching cost here is attribute-dependent (skill gap), testing the CB algorithm when Assumption~\ref{ass:key} is only an approximation.

We compare the CB algorithm against two benchmarks:
\begin{itemize}
    \item Bubble Algorithm \citep{Gan2023}: The industry-standard heuristic in gaming.
    Each player's search radius expands linearly over time; a match occurs when two radii overlap.
    This policy implicitly trades off skill gap against waiting time at the individual level.
    \item Threshold Policy: Matches the closest skill pairs whenever the queue size reaches a fixed threshold, exploiting market thickness to improve match quality.
\end{itemize}
For the CB algorithm, we track the total accumulated waiting time of all queued players as the waiting cost, and use the weighted minimum achievable skill gap as the matching cost.
A match is triggered when the accumulated realized waiting time of queued players reaches $\gamma/\alpha$ times the current minimum skill gap.

We simulate 100 player arrivals per episode, averaged over 1,000 episodes.
The same random sample paths are used across policies for each value of $\gamma$.
All algorithm parameters (expansion rate for Bubble, queue threshold for Threshold, and the empirical scaling factor $\alpha$ for CB) are optimized via grid search for each value of $\gamma$.
This empirical $\alpha$ is an implementation-specific calibration and is distinct from the robust worst-case choice $\alpha=1$ in Theorem~\ref{thm:main}.
Appendix~\ref{app:gaming_experiment} reports the parameter grids, cost decomposition, and additional calibration diagnostics for the gaming experiment.

Table~\ref{tab:compare} reports the total system costs.
The CB algorithm achieves the lowest mean cost for every reported value of $\gamma$, outperforming both the industry-standard Bubble algorithm and the fixed Threshold policy.
Relative to Bubble, CB reduces total cost by 3.66\% to 8.23\%.
Relative to Threshold, CB reduces total cost by 1.31\% to 18.41\%.
These improvements are systematic across the reported cost weights, supporting the superiority of the cost-balancing timing rule in the calibrated gaming environment.

These results illustrate two key advantages of the cost-balancing approach.
First, unlike Bubble which makes local decisions based on individual player waiting times, CB aggregates system-wide waiting costs, enabling timing decisions based on aggregate system congestion.
During congestion, accumulated waiting costs grow rapidly, triggering faster matches to clear backlogs; during lulls, the algorithm naturally waits longer to improve match quality.
Second, unlike the Threshold policy which commits to a fixed queue size, CB adapts its effective threshold based on the realized cost ratio, responding flexibly to non-stationary demand patterns.
This adaptivity explains why CB outperforms the fixed Threshold policy throughout the reported range: it can wait when additional market thickness is valuable, but it also accelerates matching when realized waiting pressure becomes high.
This suggests that the cost-balancing principle can be useful in continuous-attribute settings with non-stationary dynamics, even though the formal competitive guarantee is proved for the queue-length model.

\begin{table}[htbp!]
\centering
\begin{tabular}{ccccc}
\toprule
$\gamma$ & Bubble & CB & Threshold & Improvement over Bubble (\%) \\
\midrule
0.1 & 546.37 $\pm$ 3.72 & 505.93 $\pm$ 3.65 & 552.16 $\pm$ 3.61 & 7.40 \\
0.2 & 809.09 $\pm$ 5.99 & 742.53 $\pm$ 5.37 & 910.04 $\pm$ 7.80 & 8.23 \\
0.3 & 1021.70 $\pm$ 7.60 & 952.20 $\pm$ 7.45 & 1056.54 $\pm$ 8.72 & 6.80 \\
0.4 & 1202.25 $\pm$ 9.08 & 1121.74 $\pm$ 9.09 & 1196.39 $\pm$ 9.78 & 6.70 \\
0.5 & 1376.81 $\pm$ 11.67 & 1297.21 $\pm$ 10.93 & 1339.09 $\pm$ 10.87 & 5.78 \\
0.6 & 1502.91 $\pm$ 12.49 & 1443.53 $\pm$ 12.72 & 1473.94 $\pm$ 12.16 & 3.95 \\
0.7 & 1650.70 $\pm$ 13.10 & 1590.27 $\pm$ 13.81 & 1615.06 $\pm$ 13.14 & 3.66 \\
0.8 & 1799.49 $\pm$ 14.97 & 1732.27 $\pm$ 16.36 & 1755.35 $\pm$ 14.71 & 3.74 \\
0.9 & 1927.02 $\pm$ 15.88 & 1854.94 $\pm$ 18.46 & 1890.97 $\pm$ 16.85 & 3.74 \\
\bottomrule
\end{tabular}
\caption{\label{tab:compare}Comparison of the CB algorithm and two benchmark policies. Entries report mean total cost $\pm$ half-width of the 95\% confidence interval over 1,000 episodes.}
\end{table}

\subsection{On-Demand Last-Mile Delivery}
\label{sec:delivery}

In this subsection, we examine an on-demand last-mile delivery problem using real-world data from a food delivery platform in Shanghai, China.
The data set, adapted from \cite{LHM2021}, contains order and delivery records over a two-month period in 2015.
In this setting, drivers collect meal boxes from a central depot, deliver all assigned orders, and return to the depot before receiving new assignments.
Once dispatched, drivers cannot be rerouted to additional orders.
The service provider faces a fundamental trade-off: waiting longer allows more orders and drivers to accumulate, enabling better matching and routing efficiency, but increases customer waiting time and risks violating delivery commitments.

\subsubsection{Problem Setting}

The delivery operation involves three sequential decisions:
\begin{enumerate}
    \item Dispatch timing: Upon receiving new orders from customer locations $\mathcal{I}$ at time $t$, the service provider decides whether to dispatch the current pool of orders and drivers or to wait.
    If the decision is to wait, the system holds for a unit time interval $\Delta t$ and reevaluates.
    \item Order assignment: When a dispatch is triggered, the provider assigns a subset of locations $\mathcal{I}_k$ to each driver $k$ from the available driver pool $\mathcal{K}$.
    Each driver has a capacity constraint limiting the number of orders they can carry.
    \item Route planning: Each driver $k$ plans a delivery route that starts and ends at the central depot, visiting all assigned locations $\mathcal{I}_k$.
    The route determines the total travel time $l_k(\mathcal{I}_k)$.
\end{enumerate}

The service provider's objective is to minimize the total delivery delay across all dispatches.
Let $\tilde{t}_i$ denote the realized service time recorded for customer location $i$ (time spent at the location for handoff), and let $l_k(\mathcal{I}_k)$ denote the total travel time of driver $k$ given the assigned locations.
We treat these service times as observed inputs in the delivery experiment, since the empirical exercise focuses on dispatch timing, assignment, and routing rather than modeling service-time uncertainty.
The delivery delay for driver $k$ is the amount by which the realized total delivery time exceeds the commitment time $T_c$.
Summing over all drivers, the total delay is:
\begin{equation}
\label{eq:performance measure}
    H =\sum_{k \in \mathcal{K}}\left(\sum_{i \in \mathcal{I}_k} \tilde{t}_i+l_k\left(\mathcal{I}_k\right)-T_c\right)^{+},
\end{equation}
where $(x)^+ = \max(x, 0)$ captures that only positive delays count.
The commitment time $T_c$ is adjusted as $T_c = T_0 - \text{(waiting time)}$, where $T_0$ is the initial commitment and the waiting time accounts for how long orders have been held before dispatch.
This adjustment ensures that customers who wait longer before dispatch receive correspondingly tighter delivery windows.

The order assignment problem (Stage 2) is formulated as a vehicle routing problem (VRP) with capacity constraints, where the objective is to minimize realized total delay subject to each driver's capacity limit.
The VRP formulation, constraint descriptions, and solution methodology are provided in Appendix~\ref{app:delivery_vrp}.

This delivery setting is richer than the complete-tuple and fixed-consumption models analyzed in Sections~\ref{sec:model}--\ref{sec:extensions}.
The current pool contains pending orders and available drivers; when a dispatch is triggered, the platform selects capacity-feasible driver-route bundles through the VRP.
We therefore use the delivery study as an empirical test of the cost-balancing principle in a one-to-many dispatch environment, not as an application of the formal competitive-ratio guarantee.

\subsubsection{Experimental Setup}

The data set contains 17,645 customer orders and 3,411 driver arrivals across 839 customer locations.
We focus on lunch peak periods with batches containing at least 20 customer-order records, yielding 179 batches with an average of 86 records per batch, where each batch is 15 minutes long.
For each batch, we randomly sample order and driver arrival times within the batch window to simulate arrival variability, rather than service-time uncertainty, repeating this five times to obtain 895 arrival-varying batch samples.
We use an 80/20 calibration-evaluation split, with 716 batches used for parameter calibration and 179 held-out batches used for policy evaluation.
Key parameters are set as follows: initial commitment time $T_0 = 35$ minutes, driver capacity $C = 30$ meal boxes, and decision granularity $\Delta t = 1$ minute.

We implement the CB algorithm in this delivery context using the same operational principle: trigger a dispatch when accumulated waiting pressure reaches a calibrated proportion of the current clearing burden.
Because the delivery environment has one-to-many assignment and routing constraints, these two quantities are empirical delivery-specific proxies rather than literal instantiations of the queue-length waiting and matching costs in the theoretical model.
At each minute, the service provider evaluates whether to dispatch based on the balance between these two empirical cost components.
Let $\omega$ denote the elapsed holding time since the last dispatch.
The waiting cost is defined as $W=|\mathcal K|\omega$, the elapsed holding time multiplied by the number of available drivers.
This directly captures the reduction in remaining delivery commitment time: as waiting accumulates, the adjusted commitment time $T_c=T_0-\omega$ shrinks, increasing the risk of delay.
The matching cost is defined as the average delivery time per driver, which captures the delivery burden of the VRP clearing plan.
As more orders and drivers accumulate, denser pools can improve routing efficiency and load balancing, reducing the average clearing cost.
Appendix~\ref{app:matching_cost_monotonicity} provides empirical support for this average-cost economies pattern.
A dispatch is triggered when the current matching cost falls below $\alpha$ times the accumulated waiting cost, i.e., when $M \le \alpha W$, equivalently $M \le \alpha |\mathcal K|\omega$.
The empirical $\alpha$ is calibrated on the calibration set and absorbs the unit normalization between the waiting cost $W$ and the matching cost $M$.
The pseudocode for this implementation is provided in Appendix~\ref{app:delivery_algorithm}.

We compare the CB algorithm against four benchmark policies: (1) Greedy, which dispatches immediately when both pools are non-empty; (2) Time Threshold, which dispatches after a fixed waiting time; (3) Quantity Threshold, which dispatches when the smaller pool exceeds a fixed size; and (4) Cost-Based Z-Threshold \citep{GG2024}, which dispatches when the sum of accumulated waiting and delivery costs exceeds a threshold $Z$. We also include the Practice policy, which reflects the actual batching strategy used by the platform. 

For travel time estimation, we consider two approaches. 
The traveling salesman problem (TSP)-based approach computes optimal routes using a calibrated travel time matrix. 
The Travel Time Predictor Model, trained on actual delivery data following \cite{LHM2021}, provides more realistic estimates by capturing driver behavior that deviates from planned routes.
The rationale of adopting these two estimation approaches is discussed in Appendix~\ref{app:travel_time_discussion}.

\subsubsection{Results and Discussion}

We evaluate policy performance on the held-out evaluation set (179 batches).
To ensure a fair comparison, all policy parameters are calibrated via grid search on the same 716-batch calibration set, and all reported results are computed on the same held-out evaluation set.
This common calibration-evaluation protocol gives each benchmark a comparable tuning opportunity, avoids policy-specific cherry-picking, and helps ensure that performance differences reflect policy design rather than uneven tuning or different evaluation instances.
Details of the parameter calibration procedure are provided in Appendix~\ref{app:delivery_calibration}.

Table~\ref{tab:avg delay} reports the average delay per driver in seconds under each policy, evaluated using both the TSP-based travel time estimates and the Predictor Model estimates.
Entries report the sample mean together with the half-width of a 95\% confidence interval.
\begin{table}[ht]
    \centering
    \begin{tabular}{lcc}
    \toprule
        Policy & TSP & Predictor Model \\
    \midrule
        Greedy & 76.76 $\pm$ 8.21 & 283.52 $\pm$ 19.69 \\
        Time Threshold & 13.84 $\pm$ 3.18 & 176.23 $\pm$ 12.71 \\
        Quantity Threshold & 4.51 $\pm$ 1.26 & 122.83 $\pm$ 10.20 \\
        Cost-Based Z-Threshold & 31.84 $\pm$ 5.02 & 139.74 $\pm$ 12.75 \\
        \textbf{CB} & \textbf{4.45 $\pm$ 1.25} & \textbf{105.08 $\pm$ 9.46} \\
    \midrule
        Practice & \multicolumn{2}{c}{243.11} \\
    \bottomrule
    \end{tabular}
    \caption{Average delay comparison of the CB algorithm with benchmark policies. Entries report mean $\pm$ half-width of the 95\% confidence interval.}
    \label{tab:avg delay}
\end{table}

Several observations emerge from Table~\ref{tab:avg delay}.
First, the CB algorithm achieves the lowest mean average delay under both travel time estimation methods.
Under the Predictor Model, which more accurately captures real-world driver behavior, CB reduces average delay by 14.5\% compared to the best fixed-rule benchmark (Quantity Threshold), by 24.8\% compared to the Cost-Based Z-Threshold policy, and by 56.8\% compared to the Practice policy actually implemented by the platform.
The confidence intervals also show that the TSP results place CB and Quantity Threshold in a close range, whereas the Predictor Model results put CB at the favorable end of the calibrated policy set.

Second, the ranking of policies differs substantially between TSP and Predictor Model estimates.
Under TSP, the Quantity Threshold and CB policies perform similarly well in delay.
Under the Predictor Model, the separation between CB and the benchmark policies becomes more pronounced.
This discrepancy arises because TSP assumes drivers follow optimal routes, while in reality drivers deviate based on route familiarity, traffic conditions, and personal preferences.
The Predictor Model captures these behavioral patterns, making it a more reliable basis for policy evaluation.
The fact that CB maintains its advantage under the more realistic Predictor Model suggests that the cost-balancing timing rule remains effective when route behavior is estimated from actual delivery data.

Table~\ref{tab:avg delivery time} reports average delivery time per order in seconds, which includes both the waiting time before dispatch and the actual delivery time after dispatch.
\begin{table}[ht]
    \centering
    \begin{tabular}{lcc}
    \toprule
        Policy & TSP & Predictor Model \\
    \midrule
        Greedy & 359.26 $\pm$ 3.19 & 361.04 $\pm$ 2.15 \\
        Time Threshold & 402.41 $\pm$ 2.32 & 361.35 $\pm$ 1.69 \\
        Quantity Threshold & 365.62 $\pm$ 2.79 & 337.14 $\pm$ 1.85 \\
        Cost-Based Z-Threshold & 438.35 $\pm$ 2.57 & 403.17 $\pm$ 2.05 \\
        CB & 367.25 $\pm$ 2.54 & \textbf{331.41 $\pm$ 1.72} \\
    \midrule
        Practice & \multicolumn{2}{c}{476.79} \\
    \bottomrule
    \end{tabular}
    \caption{Average delivery time comparison of the CB algorithm with benchmark policies. Entries report mean $\pm$ half-width of the 95\% confidence interval. }
    \label{tab:avg delivery time}
\end{table}

The CB algorithm achieves the lowest average delivery time under the Predictor Model (331.41 seconds) and remains competitive under TSP (367.25 seconds).
Relative to the Predictor Model estimates, CB reduces average delivery time by 1.7\% compared to Quantity Threshold, by 17.8\% compared to Cost-Based Z-Threshold, and by 30.5\% compared to Practice.
This result complements the delay analysis: under the behavioral travel-time estimates, CB reduces violations of delivery commitments while also improving overall delivery efficiency.
The improvement stems from CB's ability to adaptively balance batch sizes.
When order and driver pools are thin, CB waits to accumulate more agents, capturing routing efficiency gains.
When pools become thick, CB dispatches promptly to avoid excessive pre-dispatch waiting.
This dynamic adjustment helps avoid batch-size choices that either sacrifice routing efficiency when batches are too small or create excessive delivery burden when batches are too large.

\begin{figure}[ht!]
    \centering
    \includegraphics[width=0.9\textwidth]{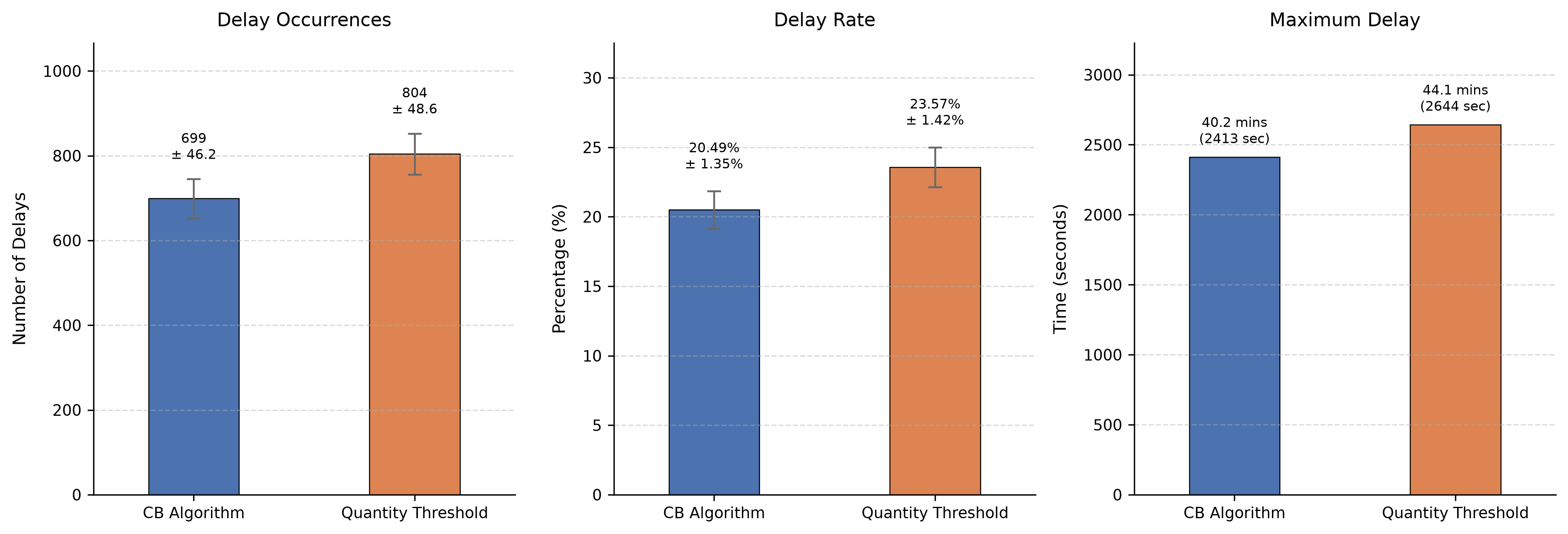}
    \caption{Performance Comparison: CB Algorithm vs Quantity Threshold}
    \label{fig:CBvsQT}
\end{figure}

Among the fixed-rule benchmarks, the Quantity Threshold policy performs best on both metrics.
This makes it the most relevant baseline for understanding CB's contribution.
Figure~\ref{fig:CBvsQT} provides a detailed comparison between CB and Quantity Threshold across multiple performance dimensions, estimated using the Travel Time Predictor Model.

The comparison reveals that CB's advantages extend beyond average performance to tail-risk metrics. 
The CB algorithm achieves a 13.1\% reduction in delay occurrences (the number of driver-routes with positive delay) and lowers the overall delay rate from 23.6\% to 20.5\%, a 3.1 percentage-point decrease.
In addition, the maximum delay under CB is 40.2 minutes, compared to 44.1 minutes under Quantity Threshold, representing an 8.8\% improvement in worst-case performance. 
These tail-risk improvements are particularly valuable in delivery operations, where extreme delays damage customer satisfaction and platform reputation disproportionately.

\begin{figure}[ht!]
    \centering
    \includegraphics[width=0.9\textwidth]{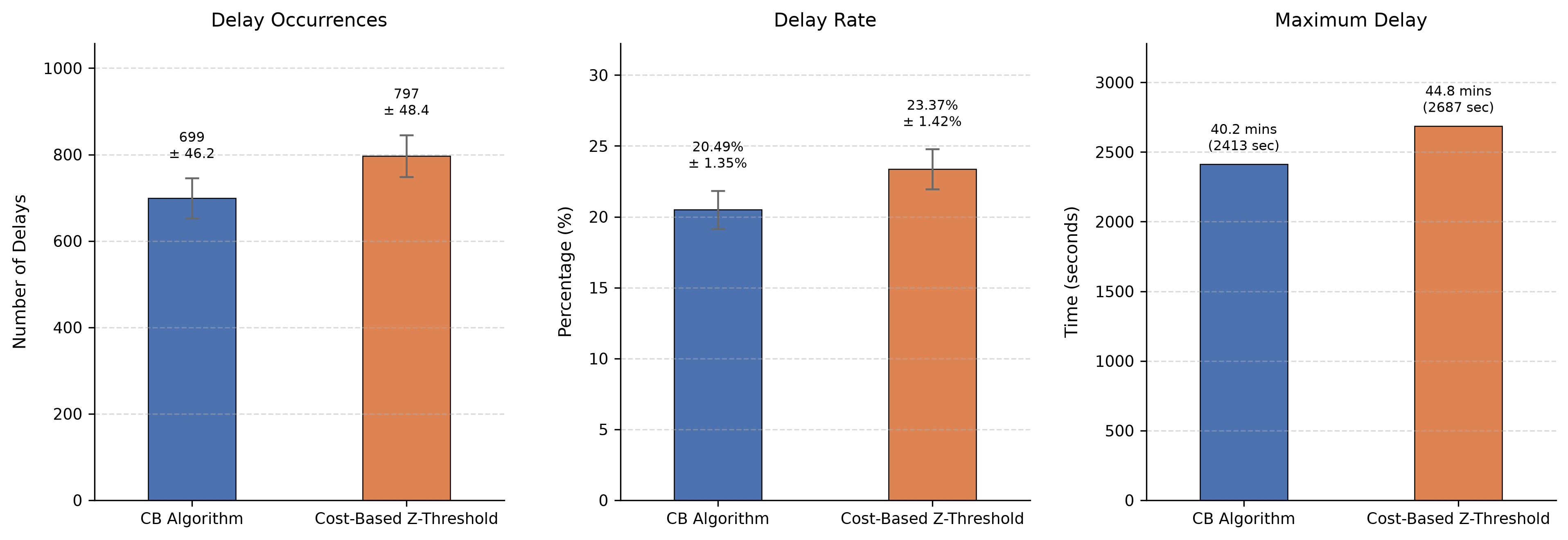}
    \caption{Performance Comparison: CB Algorithm vs Cost-Based Z-Threshold}
    \label{fig:CBvsZT}
\end{figure}

Figure~\ref{fig:CBvsZT} compares CB with the Cost-Based Z-Threshold policy proposed by \cite{GG2024}, estimated using the Travel Time Predictor Model.
This comparison is particularly instructive because both policies use cost information, but in fundamentally different ways.

The CB algorithm achieves a 12.3\% reduction in delay occurrences, measured as driver-routes with positive delay, compared to Z-Threshold; equivalently, it lowers the delay rate from 23.4\% to 20.5\%, a 2.9 percentage-point decrease. 
The key difference lies in how the two policies aggregate cost information.
The Z-Threshold policy dispatches when the sum of accumulated waiting cost and current matching cost exceeds a fixed threshold $Z$.
This additive approach treats waiting and matching costs as substitutes: a high matching cost can trigger dispatch even when waiting cost is low, and vice versa.
In contrast, CB uses a ratio-based trigger, dispatching when matching cost falls below the calibrated value of $\alpha$ times accumulated waiting cost.
This multiplicative approach treats the two costs as complements that must be balanced.

In these data, the ratio-based approach is advantageous because it responds directly to the evolving cost structure.
Consider two scenarios.
In the first scenario, matching cost is initially high because the driver pool is thin.
Under Z-Threshold, the high matching cost may push the sum above $Z$, triggering dispatch before sufficient batching efficiency is achieved.
Under CB, the ratio of matching cost to accumulated waiting cost remains high because waiting cost has not yet accumulated, so CB waits.
In the second scenario, matching cost drops quickly because drivers arrive rapidly.
Under Z-Threshold, the falling matching cost means the sum may remain below $Z$ for an extended period, which can delay dispatch even after the pool has become efficient.
Under CB, as matching cost falls, the ratio of matching cost to accumulated waiting cost decreases, eventually triggering dispatch.
In both scenarios, CB's ratio-based mechanism adapts more appropriately to the system dynamics.

\section{Conclusion}
\label{sec:conclusion}

Dynamic matching platforms face a central timing problem: matching immediately reduces waiting, while waiting can create thicker pools and improve matching efficiency. This paper develops the cost-balancing principle as a simple and robust way to manage this trade-off. The principle triggers matching when accumulated waiting cost reaches a calibrated proportion of the current matching cost. Building on this idea, we formulate a queue-length model for multi-sided matching, motivate the rule through a steady-state fluid approximation, and introduce the Cost-Balancing (CB) algorithm, which uses realized system costs rather than a fitted arrival-distribution model.

Our main contribution is to show that this simple rule has strong worst-case performance guarantees. For finite balanced complete-tuple instances, $CB_\alpha$ achieves competitive ratio $1+\max\{\alpha,1/\alpha\}$, and the robust calibration $\alpha=1$ gives a 2-competitive policy. This guarantee is distribution-free and contrasts with standard greedy and fixed-threshold policies, whose competitive ratios can be unbounded. We further show that the factor 2 is tight over the monotone complete-tuple cost class, so no deterministic online algorithm can uniformly improve on this bound in the base model. The guarantee also extends to any single fixed heterogeneous consumption requirement.

The practical value of CB lies in its state-adjusting timing rule. Unlike fixed time windows or queue-length thresholds, whose performance can deteriorate when demand conditions shift, CB changes its effective threshold with realized waiting costs, matching costs, and pool sizes. Implementation requires tracking accumulated waiting cost and evaluating the current matching or clearing cost, which can often be supplied by existing dispatch, batching, or routing modules. The numerical studies illustrate this operational value in two settings that go beyond the exact theoretical model. In the gaming experiment, CB achieves the lowest mean cost across all reported cost weights. In the food delivery experiment, CB achieves the lowest average delay under both travel-time estimates and, under the behavioral Predictor Model, also gives the lowest average delivery time. These results suggest that the cost-balancing principle captures a useful operational mechanism even in richer matching environments.

Several limitations point to future research. The current competitive guarantee assumes agents remain matchable until they are cleared, and it does not cover physical abandonment in which waiting agents leave the matching pool before being matched. Once agents can leave, different policies may face different live queues after the same arrival history, so the rank-based comparison between the $k$-th match under $CB_\alpha$ and the $k$-th match under $OPT$ no longer applies directly. The guarantee also does not cover variable clearing actions, compatibility networks, geographic restrictions, one-to-many routing, or heterogeneous preferences in full generality.
Future work could study learning-based tuning of the cost-balancing parameter, extend the guarantee to structured matching markets with richer feasibility constraints, and develop versions of cost balancing that explicitly account for abandonment before matching.

\bibliography{library}

\newpage

\clearpage
\hypersetup{pageanchor=false}
\setcounter{page}{1}
\renewcommand\thepage{\arabic{page}}
\pagestyle{headings}
\normalsize
\setcounter{nowappendix}{1}
\setcounter{section}{0}
\setcounter{subsection}{0}
\setcounter{subsubsection}{0}
\renewcommand\thesection{\Alph{section}}
\renewcommand\thesubsection{\thesection.\arabic{subsection}}
\renewcommand\thesubsubsection{\thesubsection.\arabic{subsubsection}}
\renewcommand\theparagraph{\thesubsubsection.\arabic{paragraph}}
\renewcommand\thesubparagraph{\theparagraph.\arabic{subparagraph}}

\ECHead{Online Supplement to ``When to Match: A Cost-Balancing Principle for Dynamic Markets''}

\section{Omitted Details for the Gaming Experiment}
\label{app:gaming_experiment}

This section reports additional diagnostics for the calibrated gaming experiment in Section~\ref{sec:video-game}.
Table~\ref{tab:game_decomposition} decomposes total cost into waiting and skill-gap components at the selected parameter values.
Table~\ref{tab:game_grid} reports the parameter grids used to calibrate the benchmark policies and CB.
Table~\ref{tab:game_alpha_sensitivity} reports CB's selected calibration and neighboring grid values.
Together, these diagnostics show that CB's advantage over both benchmarks is systematic across the reported gaming scenarios and is not driven by fragile parameter tuning.

\begin{table}[htbp!]
\centering
\scriptsize
\setlength{\tabcolsep}{4pt}
\begin{tabular}{cccccc}
\toprule
$\gamma$ & Policy & Waiting cost & Raw skill gap & Weighted skill gap & Total cost \\
\midrule
\multirow{3}{*}{0.1} & Bubble & 205.95 & 3404.26 & 340.43 & 546.37 \\
 & CB & 229.64 & 2762.90 & 276.29 & 505.93 \\
 & Threshold & 130.34 & 4218.23 & 421.82 & 552.16 \\
\multirow{3}{*}{0.2} & Bubble & 368.48 & 2203.04 & 440.61 & 809.09 \\
 & CB & 316.16 & 2131.88 & 426.38 & 742.53 \\
 & Threshold & 630.48 & 1397.78 & 279.56 & 910.04 \\
\multirow{3}{*}{0.3} & Bubble & 421.42 & 2000.94 & 600.28 & 1021.70 \\
 & CB & 377.72 & 1914.92 & 574.48 & 952.20 \\
 & Threshold & 633.40 & 1410.46 & 423.14 & 1056.54 \\
\multirow{3}{*}{0.4} & Bubble & 504.14 & 1745.28 & 698.11 & 1202.25 \\
 & CB & 423.28 & 1746.15 & 698.46 & 1121.74 \\
 & Threshold & 637.65 & 1396.84 & 558.74 & 1196.39 \\
\multirow{3}{*}{0.5} & Bubble & 665.61 & 1422.39 & 711.20 & 1376.81 \\
 & CB & 470.16 & 1654.10 & 827.05 & 1297.21 \\
 & Threshold & 644.00 & 1390.18 & 695.09 & 1339.09 \\
\multirow{3}{*}{0.6} & Bubble & 658.69 & 1407.04 & 844.23 & 1502.91 \\
 & CB & 535.66 & 1513.11 & 907.87 & 1443.53 \\
 & Threshold & 634.17 & 1399.61 & 839.77 & 1473.94 \\
\multirow{3}{*}{0.7} & Bubble & 661.62 & 1412.97 & 989.08 & 1650.70 \\
 & CB & 540.40 & 1499.81 & 1049.86 & 1590.27 \\
 & Threshold & 637.62 & 1396.34 & 977.44 & 1615.06 \\
\multirow{3}{*}{0.8} & Bubble & 663.81 & 1419.60 & 1135.68 & 1799.49 \\
 & CB & 576.90 & 1444.21 & 1155.37 & 1732.27 \\
 & Threshold & 639.55 & 1394.74 & 1115.79 & 1755.35 \\
\multirow{3}{*}{0.9} & Bubble & 660.44 & 1407.31 & 1266.58 & 1927.02 \\
 & CB & 641.80 & 1347.94 & 1213.14 & 1854.94 \\
 & Threshold & 638.03 & 1392.16 & 1252.94 & 1890.97 \\
\bottomrule
\end{tabular}
\caption{\label{tab:game_decomposition}Cost decomposition for the calibrated gaming experiment at selected policy parameters.}
\end{table}

Table~\ref{tab:game_decomposition} shows that CB's performance advantage does not come from mechanically minimizing only one component of the objective.
Bubble often incurs high skill-gap costs because its timing is governed by local expanding radii, while Threshold often reduces skill gaps by waiting for larger pools but pays substantial waiting costs.
CB balances these two forces directly: it waits when additional market thickness is valuable, but triggers matching once realized waiting pressure justifies the current minimum skill gap.
This cost-balance mechanism yields the lowest total cost across the reported gaming scenarios.

\begin{table}[htbp!]
\centering
\scriptsize
\begin{tabular}{lccp{0.42\textwidth}}
\toprule
Parameter & Search Range & Step & Selected values for $\gamma=0.1,\ldots,0.9$ \\
\midrule
Bubble expansion rate & 1--39 & 2 & $33,9,7,5,3,3,3,3,3$ \\
CB $1/\alpha$ & 0.05--1.95 & 0.10 & $0.65,0.75,0.75,0.75,0.75,0.85,0.75,0.75,0.85$ \\
Threshold queue size & 2--40 & 2 & $2,4,4,4,4,4,4,4,4$ \\
\bottomrule
\end{tabular}
\caption{\label{tab:game_grid}Parameter grids and selected values in the calibrated gaming experiment.}
\end{table}

Table~\ref{tab:game_grid} documents the full grid search used for each policy.
The three policies are tuned over broad, evenly spaced grids, reported by search range and step size.
The selected values are reported in the order $\gamma=0.1,0.2,\ldots,0.9$.
This calibration procedure gives each benchmark a comparable opportunity to adapt to the waiting-versus-skill-gap trade-off represented by $\gamma$.

\begin{table}[htbp!]
\centering
\scriptsize
\begin{tabular}{ccccc}
\toprule
$\gamma$ & Selected $1/\alpha$ & Lower-neighbor cost & Selected CB cost & Upper-neighbor cost \\
\midrule
0.1 & 0.65 & 506.05 & 505.93 & 506.50 \\
0.2 & 0.75 & 742.54 & 742.53 & 743.36 \\
0.3 & 0.75 & 955.70 & 952.20 & 953.93 \\
0.4 & 0.75 & 1128.57 & 1121.74 & 1122.45 \\
0.5 & 0.75 & 1300.15 & 1297.21 & 1300.15 \\
0.6 & 0.85 & 1450.47 & 1443.53 & 1446.27 \\
0.7 & 0.75 & 1597.48 & 1590.27 & 1590.53 \\
0.8 & 0.75 & 1737.49 & 1732.27 & 1737.85 \\
0.9 & 0.85 & 1865.73 & 1854.94 & 1857.05 \\
\bottomrule
\end{tabular}
\caption{\label{tab:game_alpha_sensitivity}CB $\alpha$ sensitivity in the calibrated gaming experiment. Neighboring costs use adjacent grid values around the selected $\alpha$.}
\end{table}

Table~\ref{tab:game_alpha_sensitivity} shows that the selected CB calibrations are concentrated between 0.65 and 0.85.
Moreover, neighboring grid values usually produce very similar total costs.
This pattern indicates that the gaming results are not driven by a fragile tuning choice; rather, CB performs well across a stable range of balancing parameters.

\section{Omitted Details for the Delivery Experiment}
\label{app:delivery_experiment}

This section provides additional implementation details for the delivery experiment in Section~\ref{sec:delivery}.
We formulate the underlying vehicle routing problem in Section~\ref{app:delivery_vrp} and describe its computational implementation in Section~\ref{app:vrp_implementation}, including the initial solution heuristic.
We then explain how the Cost-Balancing algorithm is implemented in the dispatching context in Section~\ref{app:delivery_algorithm} and discuss the travel time estimation methods used in the numerical study in Section~\ref{app:travel_time_discussion}.
Finally, we present the parameter calibration procedure in Section~\ref{app:delivery_calibration} and provide empirical support for the average-cost economies in Section~\ref{app:matching_cost_monotonicity}.

\subsection{VRP Formulation}
\label{app:delivery_vrp}

The order assignment problem in Section~\ref{sec:delivery} is formulated as a vehicle routing problem (VRP) with capacity constraints.
This subsection provides a detailed description of the formulation, including the meaning of each component.

Let $\mathcal{I}$ denote the set of customer locations to be served in a given dispatch batch, and let $\mathcal{K}$ denote the set of available drivers at the time of dispatch.
Each location $i \in \mathcal{I}$ has an associated order quantity $q_i$, representing the number of meal boxes to be delivered to that location.
Each driver $k \in \mathcal{K}$ has a capacity limit $C$, representing the maximum number of meal boxes the driver can carry in a single trip.

The decision variables are binary assignment variables $y_{ik} \in \{0, 1\}$ for each location-driver pair $(i, k) \in \mathcal{I} \times \mathcal{K}$.
Specifically, $y_{ik} = 1$ if location $i$ is assigned to driver $k$, and $y_{ik} = 0$ otherwise.
The vector $\mathbf{y}_k = \{y_{ik}\}_{i \in \mathcal{I}}$ collects all assignment decisions for driver $k$, defining the set of locations $\mathcal{I}_k = \{i \in \mathcal{I} : y_{ik} = 1\}$ that driver $k$ must visit.

The objective is to minimize the realized total delivery delay across all drivers for the current dispatch batch.
We treat service times as observed inputs from the delivery records rather than as a source of stochastic uncertainty in the VRP.
Let $\tilde{t}_i$ denote the realized service time at customer location $i$, which captures the time spent at the location for handoff.

For a given assignment, driver $k$'s total delivery time consists of two components: (1) the sum of realized service times at assigned locations, $\sum_{i \in \mathcal{I}} \tilde{t}_i y_{ik}$, and (2) the travel time $l_k$ required to visit all assigned locations and return to the depot.
The travel time $l_k = l(\mathbf{y}_k)$ is determined by the route selected for the assigned locations, starting at the depot, visiting all locations in $\mathcal{I}_k$, and returning to the depot.

Driver $k$'s delay is the positive part of the difference between realized total delivery time and the commitment time $T_c$:
$$
\left(\sum_{i \in \mathcal{I}} \tilde{t}_i y_{ik} + l_k - T_c\right)^{+}.
$$
The objective function sums this realized delay across all drivers:
$$
\sum_{k \in \mathcal{K}} \left(\sum_{i \in \mathcal{I}} \tilde{t}_i y_{ik} + l_k - T_c\right)^{+}.
$$

The formulation includes three types of constraints:

\begin{enumerate}
    \item Coverage constraint: Each location must be served by exactly one driver.
    $$\sum_{k \in \mathcal{K}} y_{ik} = 1, \quad \forall i \in \mathcal{I}.$$
    This ensures that every customer order is delivered and that no location is visited by multiple drivers.

    \item Capacity constraint: The total order quantity assigned to each driver must not exceed the driver's capacity.
    $$\sum_{i \in \mathcal{I}} q_i y_{ik} \leq C, \quad \forall k \in \mathcal{K}.$$
    This reflects the physical limitation that each driver can carry at most $C$ meal boxes.

    \item Routing constraint: The travel time for each driver is determined by the route selected over the assigned locations.
    $$l_k = l(\mathbf{y}_k), \quad \forall k \in \mathcal{K}.$$
    This constraint couples the assignment decisions with the routing decision.
    In practice, the implemented CVRP model uses deterministic travel-time inputs and routing constraints to compute $l_k$.
\end{enumerate}

Combining the objective and constraints, the VRP formulation is:
\begin{equation}
\label{eq:vrp}
\begin{aligned}
\min_{y_{ik}} \quad & \sum_{k \in \mathcal{K}}\left(\sum_{i \in \mathcal{I}} \tilde{t}_i y_{ik}+l_k-T_c\right)^{+} \\
\text{s.t.} \quad & \sum_{k \in \mathcal{K}} y_{ik}=1, \quad \forall i \in \mathcal{I}, \\
& \sum_{i \in \mathcal{I}} q_i y_{ik} \leq C, \quad \forall k \in \mathcal{K}, \\
& l_k=l\left(\mathbf{y}_k\right), \quad \forall k \in \mathcal{K}, \\
& y_{ik} \in\{0,1\}, \quad \forall i \in \mathcal{I}, k \in \mathcal{K}.
\end{aligned}
\end{equation}

This formulation is a mixed-integer optimization problem because assignment, routing, capacity, and delay decisions are coupled.
In the implementation, we solve deterministic CVRP instances with Gurobi, using linear delay variables for the positive-part term and deterministic travel-time inputs to evaluate route time.
Implementation details are provided in Section~\ref{app:vrp_implementation}.

\subsection{VRP Optimization Implementation}
\label{app:vrp_implementation}

This subsection provides details on the computational implementation of the VRP formulation described in Section~\ref{app:delivery_vrp}.

The VRP optimization was conducted using Gurobi Optimizer version 11.0.0 build v11.0.0rc2 (win64 - Windows 11+.0 (22631.2)). 
Computations were performed on an AMD Ryzen 5 5600G with Radeon Graphics processor featuring 6 physical cores and 12 logical processors, utilizing up to 12 parallel threads. 
Following the methodology established in \cite{LHM2021}, we implemented termination criteria of either 20 minutes of computation time or achievement of a MIP gap less than 0.01, which proved sufficient to obtain near-optimal solutions.

To enhance solution quality and computational efficiency, we developed a specialized initial solution construction heuristic.
This heuristic addresses the nontrivial case where the number of orders exceeds the number of available drivers, which is common in delivery operations during peak periods.

The heuristic operates in two phases.
In Phase 1, we initialize the assignment by distributing drivers to geographically dispersed locations.
For each driver, we identify and assign the order whose location is maximally distant from both the depot and all currently assigned order locations.
This strategy maximizes spatial coverage and ensures that drivers are initially positioned to serve orders across the service area.
Mathematically, for each driver $k$, we select order $i^* = \operatorname{argmax}_{i\in\mathcal{I}\setminus\mathcal{A}} \min_{j\in\mathcal{A}\cup\{\text{depot}\}} d(i,j)$, where $\mathcal{A}$ denotes the set of already assigned orders and $d(i,j)$ is the travel time between locations $i$ and $j$.
After Phase 1, each driver has exactly one order assigned, and we update each driver's current position and total travel time accordingly.

In Phase 2, we iteratively assign the remaining unassigned orders.
At each iteration, we identify the driver with the minimum total travel time and assign to this driver the nearest unassigned order that satisfies the capacity constraint.
This greedy strategy balances the workload across drivers while minimizing incremental travel time.
The process continues until all orders are assigned or no feasible assignment exists.

\begin{algorithm}[htbp]
\caption{Initial Solution Construction Heuristic}
\label{alg:heuristic}
\begin{algorithmic}[1]
\Require
    Order set $\mathcal{I}$, Driver set $\mathcal{K}$,
    Travel time function $d(i,j)$,
    Depot location,
    Capacity limit $C$,
    Order quantities $\{q_i\}_{i \in \mathcal{I}}$
\Ensure
    Initial feasible assignment $\mathcal{A}$ with driver positions and travel times
\State $\mathcal{A} \gets \emptyset$ \Comment{Initialize assigned orders}
\State Initialize driver positions: $\text{pos}(k) \gets \text{depot}$ for all $k \in \mathcal{K}$
\State Initialize driver travel times: $t_k \gets 0$ for all $k \in \mathcal{K}$
\State Initialize driver loads: $\text{load}(k) \gets 0$ for all $k \in \mathcal{K}$

\Statex \textbf{Phase 1: Initial Distribution}
\For{each driver $k \in \mathcal{K}$}
    \State Find order $i^* = \operatorname{argmax}\limits_{i\in\mathcal{I}\setminus\mathcal{A}} \min\limits_{j\in\mathcal{A}\cup\{\text{depot}\}} d(i,j)$
    \State Assign $i^*$ to driver $k$: $\mathcal{A} \gets \mathcal{A} \cup \{i^*\}$
    \State Update driver $k$'s position: $\text{pos}(k) \gets i^*$
    \State Update driver $k$'s travel time: $t_k \gets d(\text{depot}, i^*)$
    \State Update driver $k$'s load: $\text{load}(k) \gets q_{i^*}$
\EndFor

\Statex \textbf{Phase 2: Remaining Order Assignment}
\While{$\mathcal{A} \neq \mathcal{I}$}
    \State Find driver $k^* = \operatorname{argmin}\limits_{k\in\mathcal{K}} t_k$
    \State Find feasible nearest order: $i_{\text{nearest}} = \operatorname{argmin}\limits_{i\in\mathcal{I}\setminus\mathcal{A}: \text{load}(k^*) + q_i \leq C} d(\text{pos}(k^*), i)$
    \If{$i_{\text{nearest}}$ does not exist}
        \State \textbf{break} \Comment{No feasible assignment possible}
    \EndIf
    \State Assign $i_{\text{nearest}}$ to driver $k^*$: $\mathcal{A} \gets \mathcal{A} \cup \{i_{\text{nearest}}\}$
    \State Update driver $k^*$'s travel time: $t_{k^*} \gets t_{k^*} + d(\text{pos}(k^*), i_{\text{nearest}})$
    \State Update driver $k^*$'s position: $\text{pos}(k^*) \gets i_{\text{nearest}}$
    \State Update driver $k^*$'s load: $\text{load}(k^*) \gets \text{load}(k^*) + q_{i_{\text{nearest}}}$
\EndWhile

\State \Return Initial feasible assignment with driver positions and travel times
\end{algorithmic}
\end{algorithm}

The underlying principle of this heuristic is twofold: initially distributing drivers to geographically dispersed locations to maximize coverage, then systematically balancing total travel time by assigning orders to drivers with minimal incremental cost.
This approach not only establishes a performance floor for our VRP solution but also significantly accelerates convergence to near-optimal solutions.

To illustrate the effectiveness of our developed heuristic, we present a complex test case involving 25 drivers and 154 orders.
The corresponding VRP formulation contains 3,900 continuous variables, 604,500 binary variables, and 608,580 constraints, representing a computational challenge of considerable scale.

Table~\ref{tab:heuristic_performance} compares the computational performance of the standard Gurobi approach with our heuristic-based initialization.
Without our heuristic, Gurobi reaches the 20-minute time limit with a suboptimal objective value of 55,360.38, demonstrating the difficulty of this large-scale optimization problem.
In contrast, our heuristic generates an initial feasible solution with an objective value of 4,347.8 in under 1 second.
Using this solution as a starting point, Gurobi then converges to a near-optimal solution with an objective value of 624.83 in just 80 seconds, representing a 98.9\% improvement over the standard approach and a significant reduction in computational time.

\begin{table}[htbp]
\centering
\begin{tabular}{lcccc}
\toprule
Method & Initial Solution & Initial Time & Final Solution & Total Time \\
 & Objective Value & (seconds) & Objective Value & (seconds) \\
\midrule
Standard Gurobi & N/A & N/A & 55,360.38 & 1,200 \\
Our Heuristic + Gurobi & 4,347.8 & $<$ 1 & 624.83 & 80 \\
\midrule
Improvement & -- & -- & 98.9\% & 93.3\% \\
\bottomrule
\end{tabular}
\caption{Computational Performance Comparison for VRP Instance (25 drivers, 154 orders)}
\label{tab:heuristic_performance}
{\par \scriptsize Note: The standard Gurobi approach reached the 20-minute time limit without finding a good solution.
Our heuristic provides a high-quality starting point that enables rapid convergence to near-optimality.}
\end{table}

\subsection{CB Implementation for Delivery Dispatch}
\label{app:delivery_algorithm}

We implement the Cost-Balancing algorithm in the delivery dispatch setting as a one-to-many dispatch policy motivated by the same waiting-cost-versus-matching-cost trade-off as the main model.
At a dispatch epoch, the current pool $S$ consists of undelivered orders $\mathcal I$ and available drivers $\mathcal K$.
A feasible clearing plan assigns each driver to a capacity-feasible route serving a subset of orders, and the VRP in Section~\ref{app:delivery_vrp} computes the delivery burden of that clearing plan.
The implementation requires specifying how to measure waiting cost and matching cost in the delivery context, as well as how to implement the cost-balancing condition.

In the base model (Section~\ref{sec:model}), waiting cost accumulates at rate $\sum_{i} c_i X_i(t)$, where $X_i(t)$ is the queue length of type $i$.
The delivery implementation is not a literal instantiation of this complete-tuple cost.
Instead, it uses a dispatch-oriented urgency proxy that focuses on controllable holding time in the delivery operation.
Let $\omega$ denote the elapsed holding time, in minutes, since the last dispatch.
The waiting cost is defined as
$$
W = |\mathcal{K}| \cdot \omega,
$$
where $|\mathcal{K}|$ is the number of available drivers.
This formulation captures the reduction in remaining delivery commitment time: as waiting accumulates, the adjusted commitment time $T_c = T_0 - \omega$ shrinks, increasing the risk of delay.
Multiplying by the number of drivers reflects that more drivers waiting implies greater urgency to dispatch.

The matching cost captures the realized delivery burden per driver given the current order and driver pools.
We define it as the average delivery time per driver:
$$
M = \frac{1}{|\mathcal{K}|} \sum_{k \in \mathcal{K}}\left(\sum_{i \in \mathcal{I}_k} \tilde{t}_i+l_k\left(\mathcal{I}_k\right)\right),
$$
where $\tilde{t}_i$ is the service time at location $i$ and $l_k(\mathcal{I}_k)$ is the travel time for driver $k$ given assignment $\mathcal{I}_k$.
This quantity is the driver-normalized clearing burden used by the empirical dispatch trigger.
The intuition behind this formulation is twofold.
First, the average delivery time per order typically decreases as more orders accumulate, because routing efficiency improves with denser order distributions.
Second, the ratio of orders to drivers reflects traffic intensity: a higher ratio suggests more efficiency gain per driver from batching.
Thus, the matching cost tends to decrease over time as the system accumulates more orders and drivers, reflecting the average-cost economies used by the delivery trigger.

The CB algorithm triggers a dispatch when the matching cost falls below $\alpha$ times the accumulated waiting cost:
$$
M \le \alpha W
\quad\text{or, equivalently,}\quad
M \le \alpha |\mathcal{K}| \cdot \omega.
$$
This condition shows that the effective threshold adapts to the number of available drivers.
When more drivers are available, the algorithm dispatches sooner (lower effective threshold), preventing excessive idle time.
When fewer drivers are available, the algorithm waits longer to accumulate orders for better routing efficiency.

The calibrated value of $\alpha$ in this delivery implementation reflects the empirical normalization of $W$ and $M$.
It should not be confused with the robust calibration $\alpha=1$ in Theorem~\ref{thm:main}.

Algorithm~\ref{alg:CBM} presents the complete pseudocode for the delivery implementation of CB.
The algorithm operates in discrete time steps of $\Delta t = 1$ minute, evaluating the dispatch decision at each step.

\begin{algorithm}[htbp]
\caption{Cost-Balancing Algorithm for Delivery Dispatch}
\label{alg:CBM}
\begin{algorithmic}[1]
\Require
    Time horizon $T$, time step $\Delta t$,
    empirically calibrated balancing parameter $\alpha$
\Ensure
    Dispatch timestamps $\mathcal{D}$

\State $t \gets 0,\ \mathcal{D} \gets \emptyset,\ \omega \gets 0$
\State $\mathcal{I} \gets \emptyset$ \Comment{Current pool of undelivered orders}
\State $\mathcal{K} \gets \emptyset$ \Comment{Current pool of available drivers}
\While{$t < T$}
    \State $t \gets t + \Delta t$
    \State Update $\mathcal{I}$ with new orders arriving in $(t-\Delta t, t]$
    \State Update $\mathcal{K}$ with new drivers becoming available in $(t-\Delta t, t]$
    \If{$|\mathcal{K}| = 0$ or $|\mathcal{I}| = 0$}
        \State \textbf{continue} \Comment{No feasible dispatch}
    \EndIf
    \State $\omega \gets \omega + \Delta t$ \Comment{Elapsed holding time increment}
    \State $W \gets |\mathcal{K}| \cdot \omega$ \Comment{Accumulated waiting cost}
    \State Compute dispatch plan and obtain the current matching cost $M$
    \If{$M \le \alpha W$}
        \State $\mathcal{D} \gets \mathcal{D} \cup \{t\}$ \Comment{Record dispatch time}
        \State $\mathcal{I} \gets \emptyset$, $\mathcal{K} \gets \emptyset$
        \Comment{Remove dispatched orders and drivers}
        \State $\omega \gets 0$ \Comment{Reset elapsed holding time after dispatch}
    \EndIf
\EndWhile
\State \Return $\mathcal{D}$
\end{algorithmic}
\end{algorithm}

The algorithm initializes the system state: the current time $t$, the set of dispatch timestamps $\mathcal{D}$, the elapsed holding-time counter $\omega$, and the order and driver pools $\mathcal{I}$ and $\mathcal{K}$.
The main loop processes each time step.
At each step, the algorithm first updates the pools with newly arrived orders and drivers.
If either pool is empty, dispatch is infeasible and the algorithm continues to the next time step.
This convention treats periods with no feasible dispatch as outside the discretionary batching decision captured by the CB trigger.
When both pending orders and available drivers are present, the elapsed holding-time counter measures the controllable dispatch-holding time.
The accumulated waiting cost $W=|\mathcal{K}|\omega$ is then updated, and the current matching cost is computed by solving the VRP formulation.
The cost-balancing condition is then evaluated.
If the condition is satisfied, a dispatch is triggered: the current time is recorded, the pools are cleared, and the elapsed holding-time counter resets.
After processing all time steps, the algorithm returns the sequence of dispatch timestamps.

\subsection{Discussion on Travel Time Estimation Methods}
\label{app:travel_time_discussion}

In the experimental results, we observe notable discrepancies between the TSP-based and Predictor Model estimates, which warrant further discussion.
Understanding these differences is important because the choice of travel time estimation method affects policy evaluation and the conclusions drawn about algorithmic performance.

The TSP routing approach, widely used in both academic literature and industry practice, provides a theoretically optimal delivery path with minimum total travel time. 
In our implementation, travel times between customer locations are estimated using a calibrated travel time matrix derived from Baidu Map's RouteMatrix API. 
To account for the speeding effect from electric bikes commonly used by delivery drivers, we scale down the estimated travel time by a factor of 4/3, following the calibration approach in \cite{LHM2021}.

However, as documented in \cite{LHM2021}, real-world drivers often deviate from the theoretically optimal TSP routes.
Their empirical analysis shows that actual delivery travel time is consistently greater than or equal to the TSP solution, with an average difference of 3.8 minutes and an average relative difference of 21\%.
Several factors contribute to this systematic deviation:

\begin{enumerate}
    \item Practical road constraints: The TSP formulation does not account for real-world constraints such as limited left-turn flows at certain intersections, one-way streets, or restricted access zones that drivers must navigate around.
    \item Driver routing preferences: Drivers may prefer certain travel patterns over others.
    For instance, zigzagging routes are often avoided due to increased accident risk on busy streets, even if they represent the shortest path.
    \item Real-time adaptations: Drivers adjust their routes based on real-time traffic conditions, weather, and updated customer locations, leading to deviations from pre-planned sequences.
    \item Route familiarity: Experienced drivers may have preferred routes through familiar neighborhoods that differ from the theoretically optimal path but are more efficient in practice due to local knowledge.
\end{enumerate}

Given the difficulty of modeling all practical constraints and behavioral considerations in a mathematical formulation, we adopt a machine learning approach to predict travel time.
Following \cite{LHM2021}, we train a random forest model on the calibration set using features such as order quantity, route distance, and geographic characteristics (latitudinal and longitudinal differences between locations).
This predictor model learns from actual delivery data to estimate driver travel times without requiring specification of the visiting sequence.

The Predictor Model captures systematic differences between planned-route estimates and actual travel times by implicitly learning driver behavior patterns from historical data. 
This approach is particularly valuable because it reflects the actual consequences of routing decisions rather than idealized assumptions.

The choice of travel time estimation method has important implications for policy comparison. 
Under TSP-based estimation, policies that trigger dispatch with larger batches benefit from the idealized routing assumption, as the theoretical efficiency gains from batching are fully realized. 
Under the Predictor Model, the performance separation between CB and other policies widens because the more realistic estimates reveal performance differences that are masked by TSP assumptions.

Specifically, with large batch sizes and heavy traffic, TSP tends to underestimate delivery time because the theoretical minimum-distance route does not account for the practical constraints that become more binding as complexity increases.
Conversely, with small assignments, TSP may overestimate delivery time, as real-world drivers can leverage route familiarity and practical shortcuts.

The practice-implemented policy, which uses fixed 15-minute batching intervals, performs poorly.
This can be attributed to two factors: inflexible dispatch timing that ignores the current system state, and potentially unbalanced order assignments that lead to large variations in delivery times among drivers.

Since our goal is to evaluate policies under realistic conditions, we consider the Predictor Model results as the primary basis for comparison, while reporting TSP results for completeness and comparability with other studies that use theoretical routing assumptions.

\subsection{Parameter Calibration}
\label{app:delivery_calibration}

This subsection describes the parameter calibration procedure for all dispatch policies, including the search methodology, optimal parameter values, and interpretation of results.
Our goal is to give each policy the best chance to perform well so that performance differences reflect algorithmic design rather than poor tuning.

To ensure a fair comparison across policies, we calibrate all dispatch-policy parameters using grid search on the same calibration set (716 batches).
For each policy, we evaluate performance across a predefined parameter range and select the value that minimizes average delay on the calibration data.
This approach avoids policy-specific cherry-picking and ensures that each benchmark operates at its best possible configuration under the same calibration metric.

The grid search is conducted independently for each policy, using the same calibration-evaluation split and evaluation metric (average delay per driver).
After calibration, all policy comparisons are reported on the common held-out evaluation set (179 batches), so differences across policies are not driven by different tuning samples or different evaluation instances.
All experiments use the Travel Time Predictor Model for consistency, as it provides more realistic estimates of driver behavior than idealized TSP routing (see Section~\ref{app:travel_time_discussion}).

Table~\ref{tab:Parameter tuning} summarizes the search ranges, step sizes, and optimal values for each policy.

\begin{table}[ht]
    \centering
\begin{tabular}{lccc}
    \toprule
    Parameter & Search Range & Step & Optimal \\
    \midrule
    Time Threshold & 1--15 min & 1 min & 7 min \\
    Quantity Threshold & 1--10 & 1 & 3 \\
    Cost-Based Z-Threshold & 20--35 min & 1 min & 32 min \\
    $1/\alpha$ (CB) & 0.1--1.0 & 0.1 & 0.5 \\
    \bottomrule
    \end{tabular}
    \caption{Parameter Calibration Results for Delivery Experiment}
    \label{tab:Parameter tuning}
\end{table}

For CB, Table~\ref{tab:Parameter tuning} reports the grid in terms of $1/\alpha$ for readability.
Under the manuscript convention $M\le \alpha W$, with $W=|\mathcal K|\omega$ in the delivery implementation, the selected value $1/\alpha=0.5$ corresponds to $\alpha=2$.
This empirically calibrated value is favorable on the calibration set and consistent with an average-cost interpretation of the delivery implementation.
Specifically, the value aligns with an average-based approximation, where the waiting and matching costs are summarized over the calibration instances.
In this approximation, the number of drivers is replaced by the optimal quantity threshold, the waiting cost is replaced by the optimal waiting time threshold, and the matching cost is estimated using the sample average of delivery time and order quantity.
$$
\alpha \approx \frac{\text{average matching cost}}{\text{average waiting cost}} = \frac{851.21}{7 \cdot 60} \approx 2.03.
$$
Under this average-based approximation, the derived $\alpha$ value also approximates 2.
This calibration reflects the empirical units used for $W$ and $M$ in the delivery study and is distinct from the robust base-model calibration $\alpha=1$ in Theorem~\ref{thm:main}.

The initial delivery commitment time $T_0$ is set at 35 minutes, which is slightly stricter than the 45-minute setting in \cite{LHM2021}.
This stricter commitment time simulates a high-load environment that better highlights the differences in algorithmic performance, as policies must operate closer to their capacity limits.

The travel time predictor model is trained on the 716 calibration batches using a random forest algorithm with features proposed by \cite{LHM2021}.
These attributes summarize both the scale and geometry of each batch of customer locations.
For completeness, we include the feature table below.

\begin{table}[htbp]
    \centering
    \begin{tabular}{p{0.18\linewidth}p{0.75\linewidth}}
    \toprule
    Attribute & Definition \\
    \midrule
    $\bar{d}$ & Average distance between the depot and customer locations \\
    $d$       & Shortest distance between the depot and customer locations \\
    $D$       & Longest distance between the depot and customer locations \\
    $R$       & Area of the smallest rectangle covering the customer locations \\
    $R'$      & Area of the smallest rectangle covering both the depot and the customer locations \\
    $L$       & Area of the smallest lune (formed by two equal-angle sectors from the depot) covering the customer locations \\
    $a$       & Maximum latitudinal difference between a pair of customer locations \\
    $b$       & Maximum longitudinal difference between a pair of customer locations \\
    $a'$      & Maximum latitudinal difference between any two locations including the depot \\
    $b'$      & Maximum longitudinal difference between any two locations including the depot \\
    $\mathrm{cstdev\_a}$ & Standard deviation of latitudinal differences between the depot and customer locations \\
    $\mathrm{cstdev\_b}$ & Standard deviation of longitudinal differences between the depot and customer locations \\
    $s_a$     & Average latitudinal difference between pairs of customer locations \\
    $s_b$     & Average longitudinal difference between pairs of customer locations \\
    \bottomrule
    \end{tabular}
    \caption{Graph-Based Attributes Used in the Travel Time Predictor Model}
    \label{tab:predictor_features}
\end{table}

\subsection{Empirical Support for Average-Cost Economies}
\label{app:matching_cost_monotonicity}

The delivery trigger is motivated by average-cost economies: as the current pool becomes thicker, the average dispatch burden can decrease while waiting cost accumulates.
We examine this pattern empirically using 2,323 dispatch batches from our delivery experiment, sampled from the decision points of the optimally implemented CB algorithm with the manuscript-convention value $\alpha = 2$.

We first examine how matching cost varies with the number of accumulated orders.
In the delivery context, larger order pools enable more efficient routing: drivers can serve multiple nearby locations in a single trip, reducing travel time per order.
Table~\ref{tab:monotonicity_order} presents the grouped analysis by order quantity.

\begin{table}[htbp]
    \centering
    \begin{tabular}{cccc}
    \toprule
    Order Range & Mean Orders & Sample Size & Avg Delivery Time (s) \\
    \midrule
    1--12 & 6.9 & 1205 & 521.2 \\
    13--23 & 17.1 & 714 & 319.5 \\
    24--34 & 27.6 & 241 & 261.6 \\
    35--46 & 39.1 & 98 & 225.2 \\
    47--57 & 50.8 & 39 & 202.1 \\
    58--68 & 62.4 & 18 & 218.7 \\
    70--79 & 74.0 & 6 & 186.3 \\
    \bottomrule
    \end{tabular}
    \caption{Average Delivery Time by Order Quantity}
    \label{tab:monotonicity_order}
\end{table}

The data reveals a clear decreasing trend: as the number of orders increases from 1--12 to 70--79, the average delivery time decreases from 521 seconds to 186 seconds, representing a 64.3\% reduction.
The corresponding simple linear regression is
$$
\text{avg\_delivery\_time} = 560.81 - 10.00 \times \text{order\_num},
$$
indicating that each additional order reduces the average delivery time by approximately 10 seconds.
This substantial improvement reflects the routing efficiency gains from denser order distributions.

We note a slight non-monotonicity in the 58--68 order range (218.7 seconds), which is higher than the adjacent 47--57 range (202.1 seconds).
This deviation is likely due to limited sample size (only 18 observations) and the inherent variability in order locations.
Despite this local fluctuation, the overall trend is clearly decreasing, and the regression captures the dominant relationship.

Using order quantity alone does not fully capture the system state, as the driver pool also affects matching efficiency.
When more drivers are available, orders can be distributed more evenly, reducing the workload per driver.
To capture the joint effect of order and driver accumulation, we analyze the relationship between waiting time and matching cost.
As the system waits longer, both order and driver pools grow, and their interaction determines the actual matching cost.

\begin{table}[htbp]
    \centering
    \begin{tabular}{cccc}
    \toprule
    Waiting Time (min) & Sample Size & Avg Delivery Time (s) & Avg Orders \\
    \midrule
    1 & 489 & 531.9 & 7.2 \\
    2 & 608 & 435.1 & 13.6 \\
    3 & 570 & 388.9 & 18.4 \\
    4 & 366 & 344.4 & 19.0 \\
    5 & 173 & 302.7 & 19.2 \\
    6 & 68 & 267.2 & 20.0 \\
    7 & 31 & 247.6 & 18.4 \\
    8 & 14 & 224.0 & 22.1 \\
    9 & 4 & 195.5 & 20.0 \\
    \bottomrule
    \end{tabular}
    \caption{Average Delivery Time by Waiting Time}
    \label{tab:monotonicity_waiting}
\end{table}

Table~\ref{tab:monotonicity_waiting} shows that average delivery time decreases monotonically from 532 seconds at 1 minute of waiting to 195 seconds at 9 minutes, a 63.3\% reduction.
Unlike the order quantity analysis, the waiting time analysis exhibits strict monotonicity across all observed values.
The corresponding regression is
$$
\text{avg\_delivery\_time} = 554.33 - 51.13 \times \text{waiting\_time},
$$
indicating that each additional minute of waiting reduces the average delivery time by approximately 51 seconds.
This larger coefficient (compared to 10 seconds per order) reflects that waiting time captures multiple effects simultaneously: more orders, more drivers, and better order-driver matching opportunities.

\begin{figure}[htbp]
    \centering
    \includegraphics[width=0.9\textwidth]{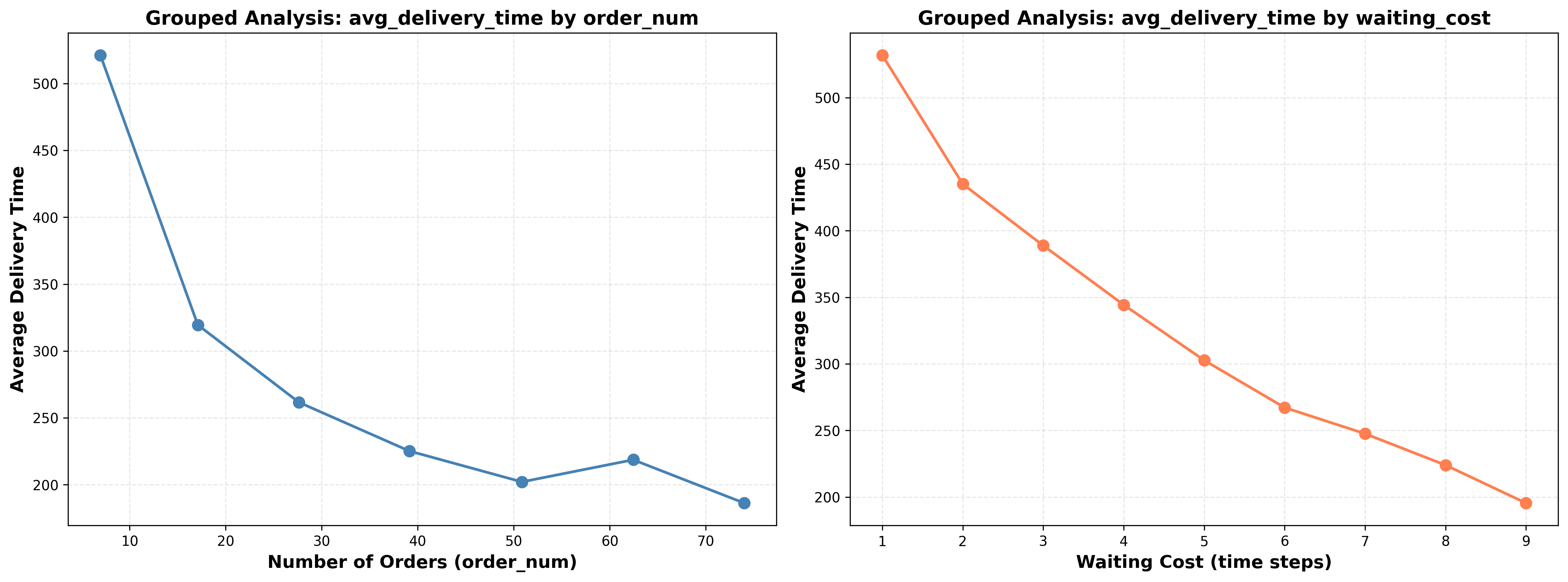}
    \caption{Empirical Support for Average-Cost Economies}
    \label{fig:monotonicity}
\end{figure}

Figure~\ref{fig:monotonicity} visualizes these relationships.
The left panel shows the decreasing trend of average delivery time with order quantity, while the right panel demonstrates the monotonic decrease with waiting time.
A notable feature of both curves is that the slope is steep at low values and gradually flattens as queue length or waiting time increases. 
This pattern reflects diminishing marginal returns: the benefit of adding one more agent or waiting one more minute is largest when the pool is small and decreases as more agents accumulate. 
In the left panel, as we can find from Table~\ref{tab:monotonicity_order}, moving from 7 to 17 orders reduces delivery time by approximately 200 seconds, while moving from 50 to 74 orders yields only about 16 seconds of improvement.
Similarly, in the right panel, the first few minutes of waiting produce the largest gains.

This empirical pattern supports the average-cost economies used to motivate the delivery implementation.
The observed curvature in Figure~\ref{fig:monotonicity} also suggests diminishing marginal returns: the first few additional orders or minutes of waiting generate the largest reductions in average delivery time.

These results provide empirical support for treating the delivery setting as a one-to-many dispatch environment with average-cost economies.
The economies-of-scale effect arises from two complementary mechanisms: (1) routing efficiency, where denser order distributions enable shorter travel distances per delivery, and (2) load balancing, where larger driver pools allow more even distribution of orders across drivers.
Together, these mechanisms create substantial benefits from waiting, which the CB algorithm exploits by delaying dispatch until the accumulated waiting cost justifies the current matching cost.
The evidence here supports the average-cost mechanism used in the delivery implementation, while the formal competitive guarantee remains tied to the complete-tuple and fixed-consumption models in the main text.

\section{Discussion on Model with Attribute Information}
\label{sec:general_model}

Our main analysis rests on Assumption~\ref{ass:key}, which posits that matching costs depend on queue lengths (market thickness) rather than the microscopic attributes of individual agents.
While real-world costs are indeed attribute-dependent, it is natural to ask what happens when full attribute information is incorporated into the online decision-making framework.
In this section, we formally define the general attribute-based model and characterize two fundamental barriers that arise: computational intractability and theoretical impossibility of bounded competitive ratios under adversarial arrivals.

\subsection{The General Attribute-Based Model}
\label{sec:attribute-model}

We now formally define the general matching model with attribute information, which extends the queue-based model in Section~\ref{sec:problem definition}.

In the attribute-based model, each arriving agent carries relevant attributes drawn from an attribute space $\mathcal{Q}$.
Specifically, the $n$-th type-$i$ agent arrives with attribute $q_{in} \in \mathcal{Q}$ drawn from a type-specific distribution $\nu_i$.
These attributes encode agent-specific information relevant to matching quality: geographic locations in ride-sharing, skill ratings in gaming, or delivery destinations in logistics.

The system state at time $t$ consists of two components:
\begin{enumerate}
    \item The queue length vector $\mathbf{X}(t) = (X_1(t), \ldots, X_N(t)) \in \mathbb{Z}_+^N$, where $X_i(t)$ is the number of unmatched type-$i$ agents.
    \item The attribute configuration $\mathbf{Q}(t) = \{q_{ij} : j \leq X_i(t), i \in \{1, \ldots, N\}\}$, collecting the attributes of all unmatched agents.
\end{enumerate}

In this general model, the matching cost depends on the specific agents selected for a match.
Let $f(q_1, \ldots, q_N)$ denote the cost of matching a tuple of agents with attributes $(q_1, \ldots, q_N)$.
When the controller decides to match at state $(\mathbf{X}, \mathbf{Q})$, it selects one agent from each type to minimize the matching cost:
$$
f^*(\mathbf{X}, \mathbf{Q}) := \min_{j_i \leq X_i, \forall i} f(q_{1,j_1}, \ldots, q_{N,j_N}).
$$
A policy $\pi$ in this model is a mapping from the full state $(\mathbf{X}, \mathbf{Q}, t)$ to a matching decision.
The monotonicity property of the matching cost function is modified accordingly for the general attribute-based model.
\begin{assumption}
    \label{ass:general_attribute-based_model}
    The matching cost incurred at time $t$ is a function $f^*(\cdot, \cdot)$, which depends on $\{X_i(t), 1\le i\le N\}$, the number of unmatched agents at time $t$ for each type, and the attributes $\{Q_i(t), 1\le i\le N\}$ of all unmatched agents.
    Furthermore, $f^*(\cdot, \cdot)$ satisfies the monotonicity property, i.e., $f^*(\mathbf{x}, \mathbf{q}_x)\le f^*(\mathbf{y}, \mathbf{q}_y)$ whenever $\mathbf{x} \succeq \mathbf{y}$ and the attributes $\mathbf{q}_y$ are a subcollection of $\mathbf{q}_x$ (i.e., the set of agents in the queues for $\mathbf{y}$ is a subset of those for $\mathbf{x}$).
\end{assumption}

The queue-based model in our main analysis (Assumption~\ref{ass:key}) can be viewed as a tractable approximation where the matching cost $f(\mathbf{X})$ depends only on queue lengths, not on specific attribute realizations.
This approximation is justified when larger pools statistically yield better matches, which holds in most practical settings.
The following results characterize the computational and information-theoretic barriers inherent to this general formulation.

\subsection{Computational Intractability}
\label{sec:computational_intractability}

We first address computational feasibility.
Even in a static setting where all agents are present simultaneously, finding the cost-minimizing match is prohibitively difficult for multi-sided markets.

\begin{proposition}
    \label{prop:nphard}
    For $N \geq 3$, computing an optimal matching in the attribute-based model is NP-hard, even in a static one-shot setting with zero waiting costs.
\end{proposition}

\begin{remark}
    The hardness result extends to any $N \geq 3$ via reduction to $N$-dimensional matching.
    For $N = 2$, the static assignment is polynomial-time solvable as minimum-weight bipartite matching \citep{Kuhn1955}, but the dynamic problem remains intractable due to the continuous attribute state space and infinite horizon.
\end{remark}

This result has important practical implications.
For platforms coordinating three or more agent types (e.g., rider, driver, and restaurant in food delivery), even computing the optimal match for a fixed set of agents is computationally prohibitive.
Dynamic optimization over an infinite horizon is even more demanding.
By abstracting attributes into aggregate queue lengths, we bypass this combinatorial explosion, focusing on the temporal trade-offs that are central to dynamic matching control.

\subsection{Theoretical Impossibility Under Adversarial Attributes}
\label{sec:impossibility}

Beyond computational hardness, there is a more fundamental barrier: the impossibility of bounded competitive ratios for deterministic online algorithms when the adversary controls agent attributes.
In our main model, adversarial arrivals are restricted to timing; attributes are abstracted into queue lengths.
Here we show that if the adversary can also choose specific agent attributes, no deterministic online algorithm can achieve bounded performance guarantees.

\begin{proposition}
    \label{prop:adversarial attribute}
    In the attribute-based model where the adversary controls both arrival times and agent attributes, no deterministic online algorithm achieves a bounded competitive ratio.
\end{proposition}

This impossibility result reveals a fundamental tension: no deterministic online algorithm can simultaneously be robust to adversarial attributes and adversarial arrival timing.
The adversary can always exploit the algorithm's ignorance of future attribute realizations to create arbitrarily bad outcomes.

These results clarify the role of the queue-based abstraction (Assumption~\ref{ass:key}).
By modeling matching costs as functions of queue lengths rather than individual attributes, the queue-based model captures the empirical regularity that larger pools yield better matches on average, without requiring the algorithm to predict specific attribute realizations.
This abstraction enables the bounded competitive ratio guarantee of Theorem~\ref{thm:main} for CB by restricting the adversary's power to arrival timing while retaining monotone queue-length matching costs.

\clearpage
\setcounter{page}{1}
\ECHead{Online Appendix to ``When to Match: A Cost-Balancing Principle for Dynamic Markets''}

\section{Main Proofs}
\label{app:proofs}

\subsection{Proof of Proposition~\ref{prop:relaxed_monotonicity}}
\label{app:proof_relaxed_monotonicity}

To prove the proposition, we first establish the following technical lemma.
\begin{lemma}
    \label{lemma:convex_supermodular_non_increasing}
    If a function $f(\mathbf{x})$ is convex, supermodular, and component-wise non-increasing, then the function $g(\mathbf{x}) = f(\mathbf{x} + \mathbf{e}_i) - f(\mathbf{x})$ is non-decreasing for all $\mathbf{x}$ and $i$.
\end{lemma}

\proof{Proof of Lemma~\ref{lemma:convex_supermodular_non_increasing}.}
    We want to show that $g(\mathbf{x}) = f(\mathbf{x} + \mathbf{e}_i) - f(\mathbf{x})$ is non-decreasing in $\mathbf{x}$, i.e., $g(\mathbf{y}) \ge g(\mathbf{x})$ for all $\mathbf{y} \ge \mathbf{x}$.
    It suffices to show that $g(\mathbf{x} + \mathbf{e}_j) \ge g(\mathbf{x})$ for all $j$.

    \textbf{Case 1: $j=i$.}
    $$g(\mathbf{x} + \mathbf{e}_i) - g(\mathbf{x}) = [f(\mathbf{x} + 2\mathbf{e}_i) - f(\mathbf{x} + \mathbf{e}_i)] - [f(\mathbf{x} + \mathbf{e}_i) - f(\mathbf{x})].$$
    From the convexity of $f$, we have $2f(\mathbf{x} + \mathbf{e}_i) \le f(\mathbf{x}) + f(\mathbf{x} + 2\mathbf{e}_i)$, which implies $f(\mathbf{x} + 2\mathbf{e}_i) - f(\mathbf{x} + \mathbf{e}_i) \ge f(\mathbf{x} + \mathbf{e}_i) - f(\mathbf{x})$.
    Thus, $g(\mathbf{x} + \mathbf{e}_i) \ge g(\mathbf{x})$.

    \textbf{Case 2: $j \ne i$.}
    $$g(\mathbf{x} + \mathbf{e}_j) - g(\mathbf{x}) = [f(\mathbf{x} + \mathbf{e}_i + \mathbf{e}_j) - f(\mathbf{x} + \mathbf{e}_j)] - [f(\mathbf{x} + \mathbf{e}_i) - f(\mathbf{x})].$$
    From the supermodularity of $f$, we have $f(\mathbf{x} + \mathbf{e}_i) + f(\mathbf{x} + \mathbf{e}_j) \le f(\mathbf{x}) + f(\mathbf{x} + \mathbf{e}_i + \mathbf{e}_j)$, which implies $f(\mathbf{x} + \mathbf{e}_i + \mathbf{e}_j) - f(\mathbf{x} + \mathbf{e}_j) \ge f(\mathbf{x} + \mathbf{e}_i) - f(\mathbf{x})$.
    Thus, $g(\mathbf{x} + \mathbf{e}_j) \ge g(\mathbf{x})$.

    Combining both cases, $g(\mathbf{x})$ is non-decreasing in $\mathbf{x}$.
\hfill\Halmos
\endproof

\vspace{1em}

We prove the finite-horizon statement by analyzing the structural properties of the value function in the uniformized dynamic program.
Fix a uniformization rate $\Lambda>\sum_{i=1}^N\lambda_i$ and, for any function $g$, define
\begin{equation}
\label{eq:uniformized_operator}
    \mathcal U g(\mathbf{x})
    =
    \left(1-\sum_{i=1}^N\frac{\lambda_i}{\Lambda}\right)g(\mathbf{x})
    +
    \sum_{i=1}^N \frac{\lambda_i}{\Lambda}g(\mathbf{x}+\mathbf e_i).
\end{equation}
Thus $\mathcal U$ is the transition operator over one uniformized decision epoch, including the dummy no-arrival event.
Let $J_k(\mathbf{x})$ denote the optimal $k$-stage cost, initialized by $J_0(\mathbf{x})=0$.
Given $J_k$, define the wait and match action values for the $(k+1)$-stage problem by
\begin{align}
\label{eq:uniformized_wait_value}
    J_{k+1}^0(\mathbf{x})
    &=
    \frac{w(\mathbf{x})}{\Lambda}
    +
    \mathcal U J_k(\mathbf{x}),\\
\label{eq:uniformized_match_value}
    J_{k+1}^1(\mathbf{x})
    &=
    f(\mathbf{x})
    +
    \frac{w(\mathbf{x}-\mathbf 1)}{\Lambda}
    +
    \mathcal U J_k(\mathbf{x}-\mathbf 1),
    \qquad \mathbf{x}\in\mathcal F .
\end{align}
If $\mathbf{x}\notin\mathcal F$, matching is infeasible and waiting is forced.
For $\mathbf{x}\in\mathcal F$,
\[
    J_{k+1}(\mathbf{x})
    =
    \min\{J_{k+1}^0(\mathbf{x}),J_{k+1}^1(\mathbf{x})\}.
\]
We select match whenever it is weakly optimal.
Equivalently, define the local net benefit of matching by
\[
    B_{k+1}(\mathbf{x})
    =
    J_{k+1}^0(\mathbf{x})-J_{k+1}^1(\mathbf{x}),
    \qquad \mathbf{x}\in\mathcal F,
\]
so matching is selected if and only if $B_{k+1}(\mathbf{x})\ge0$.
Since $J_{k+1}^1(\mathbf{x})=f(\mathbf{x})+J_{k+1}^0(\mathbf{x}-\mathbf 1)$, the diagonal difference
\[
    G_k(\mathbf{x})=J_k(\mathbf{x})-J_k(\mathbf{x}-\mathbf 1),
    \qquad \mathbf{x}\in\mathcal F,
\]
gives
\begin{equation}
\label{eq:benefit_identity}
    B_{k+1}(\mathbf{x})
    =
    \frac{\sum_{i=1}^N c_i}{\Lambda}
    -
    f(\mathbf{x})
    +
    \mathcal U G_k(\mathbf{x}).
\end{equation}
All expressions involving a match are evaluated only at feasible states; if a comparison state lies outside $\mathcal F$, the match action is unavailable and the wait action is forced.

We prove by induction that, for every $k\ge0$,
\begin{equation}
\label{eq:diagonal_bound}
    G_k(\mathbf{x}+\mathbf e_i)-G_k(\mathbf{x})
    \ge
    f(\mathbf{x}+\mathbf e_i)-f(\mathbf{x}),
    \qquad \mathbf{x}\in\mathcal F,\ i=1,\ldots,N.
\end{equation}
The base case follows from $J_0\equiv0$, because $G_0\equiv0$ and $f$ is component-wise non-increasing.
Assume \eqref{eq:diagonal_bound} holds for $k$.
Using \eqref{eq:benefit_identity}, for any $\mathbf{x}\in\mathcal F$,
\[
\begin{aligned}
    B_{k+1}(\mathbf{x}+\mathbf e_i)-B_{k+1}(\mathbf{x})
    &=
    -\bigl(f(\mathbf{x}+\mathbf e_i)-f(\mathbf{x})\bigr)
    +
    \mathcal U\bigl[G_k(\mathbf{x}+\mathbf e_i)-G_k(\mathbf{x})\bigr]  \\
    &\ge
    -\bigl(f(\mathbf{x}+\mathbf e_i)-f(\mathbf{x})\bigr)
    +
    \mathcal U\bigl[f(\mathbf{x}+\mathbf e_i)-f(\mathbf{x})\bigr] \\
    &\ge 0.
\end{aligned}
\]
The last inequality follows from Lemma~\ref{lemma:convex_supermodular_non_increasing}, because $f(\mathbf{x}+\mathbf e_i)-f(\mathbf{x})$ is non-decreasing in $\mathbf{x}$ and $\mathcal U$ is a convex combination of shifted values.
Therefore $B_{k+1}$ is non-decreasing in $\mathbf{x}$.

It remains to verify \eqref{eq:diagonal_bound} for $k+1$.
Fix $\mathbf{x}\in\mathcal F$ and $i$.
Along the chain
\[
    \mathbf{x}-\mathbf 1
    \preceq
    \mathbf{x}-\mathbf 1+\mathbf e_i
    \preceq
    \mathbf{x}
    \preceq
    \mathbf{x}+\mathbf e_i,
\]
the selected action can switch from wait to match at most once, because $B_{k+1}$ is non-decreasing and match is selected exactly when $B_{k+1}\ge0$.
Thus the following five cases cover all possibilities, with action $0$ denoting wait and action $1$ denoting match.
If a state in the table is infeasible, only cases assigning action $0$ to that state can occur.
\begin{table}[htbp]
\caption{Selected actions at neighboring states}
\label{tab:optimal_actions}
\centering
\begin{tabular}{lccccc}
\toprule
State & Case 1 & Case 2 & Case 3 & Case 4 & Case 5 \\
\midrule
$\mathbf{x}+\mathbf e_i$ & 0 & 1 & 1 & 1 & 1 \\
$\mathbf{x}$ & 0 & 0 & 1 & 1 & 1 \\
$\mathbf{x}-\mathbf 1+\mathbf e_i$ & 0 & 0 & 0 & 1 & 1 \\
$\mathbf{x}-\mathbf 1$ & 0 & 0 & 0 & 0 & 1 \\
\bottomrule
\end{tabular}
\end{table}

\begin{itemize}
    \item \textbf{Case 1:} Wait at all four states.
    Then
    \[
        G_{k+1}(\mathbf{x})
        =
        \frac{\sum_i c_i}{\Lambda}
        +
        \mathcal U G_k(\mathbf{x}),
    \]
    and hence
    \[
        G_{k+1}(\mathbf{x}+\mathbf e_i)-G_{k+1}(\mathbf{x})
        =
        \mathcal U\bigl[G_k(\mathbf{x}+\mathbf e_i)-G_k(\mathbf{x})\bigr].
    \]
    The induction hypothesis and Lemma~\ref{lemma:convex_supermodular_non_increasing} imply this quantity is at least
    $f(\mathbf{x}+\mathbf e_i)-f(\mathbf{x})$.

    \item \textbf{Case 2:} Match only at $\mathbf{x}+\mathbf e_i$.
    Since $J_{k+1}^1(\mathbf{x}+\mathbf e_i)=f(\mathbf{x}+\mathbf e_i)+J_{k+1}^0(\mathbf{x}-\mathbf 1+\mathbf e_i)$,
    \[
        G_{k+1}(\mathbf{x}+\mathbf e_i)=f(\mathbf{x}+\mathbf e_i).
    \]
    Therefore
    \[
    \begin{aligned}
        G_{k+1}(\mathbf{x}+\mathbf e_i)-G_{k+1}(\mathbf{x})
        &=
        f(\mathbf{x}+\mathbf e_i)-J_{k+1}^0(\mathbf{x})+J_{k+1}^0(\mathbf{x}-\mathbf 1)\\
        &=
        f(\mathbf{x}+\mathbf e_i)-f(\mathbf{x})
        +
        J_{k+1}^1(\mathbf{x})-J_{k+1}^0(\mathbf{x})\\
        &\ge
        f(\mathbf{x}+\mathbf e_i)-f(\mathbf{x}),
    \end{aligned}
    \]
    because wait is selected at $\mathbf{x}$.

    \item \textbf{Case 3:} Match at $\mathbf{x}+\mathbf e_i$ and $\mathbf{x}$, and wait at the two lower states.
    Then
    \[
        G_{k+1}(\mathbf{x}+\mathbf e_i)=f(\mathbf{x}+\mathbf e_i),
        \qquad
        G_{k+1}(\mathbf{x})=f(\mathbf{x}),
    \]
    so \eqref{eq:diagonal_bound} holds with equality.

    \item \textbf{Case 4:} Wait only at $\mathbf{x}-\mathbf 1$.
    Then
    \[
    \begin{aligned}
        G_{k+1}(\mathbf{x}+\mathbf e_i)-G_{k+1}(\mathbf{x})
        &=
        J_{k+1}^1(\mathbf{x}+\mathbf e_i)
        -J_{k+1}^1(\mathbf{x}-\mathbf 1+\mathbf e_i)
        -J_{k+1}^1(\mathbf{x})
        +J_{k+1}^0(\mathbf{x}-\mathbf 1)\\
        &=
        f(\mathbf{x}+\mathbf e_i)-f(\mathbf{x})
        +
        J_{k+1}^0(\mathbf{x}-\mathbf 1+\mathbf e_i)
        -J_{k+1}^1(\mathbf{x}-\mathbf 1+\mathbf e_i)\\
        &\ge
        f(\mathbf{x}+\mathbf e_i)-f(\mathbf{x}),
    \end{aligned}
    \]
    because match is selected at $\mathbf{x}-\mathbf 1+\mathbf e_i$.

    \item \textbf{Case 5:} Match at all four states.
    This case can occur only when $\mathbf{x}-\mathbf 1\in\mathcal F$.
    Direct substitution gives
    \[
        G_{k+1}(\mathbf{x})
        =
        \frac{\sum_i c_i}{\Lambda}
        +
        f(\mathbf{x})-f(\mathbf{x}-\mathbf 1)
        +
        \mathcal U G_k(\mathbf{x}-\mathbf 1).
    \]
    Hence
    \[
    \begin{aligned}
        &G_{k+1}(\mathbf{x}+\mathbf e_i)-G_{k+1}(\mathbf{x})\\
        &=
        f(\mathbf{x}+\mathbf e_i)-f(\mathbf{x})
        -f(\mathbf{x}-\mathbf 1+\mathbf e_i)+f(\mathbf{x}-\mathbf 1)\\
        &\quad
        +
        \mathcal U\bigl[
            G_k(\mathbf{x}-\mathbf 1+\mathbf e_i)-G_k(\mathbf{x}-\mathbf 1)
        \bigr]\\
        &\ge
        f(\mathbf{x}+\mathbf e_i)-f(\mathbf{x})
        -f(\mathbf{x}-\mathbf 1+\mathbf e_i)+f(\mathbf{x}-\mathbf 1)\\
        &\quad
        +
        \mathcal U\bigl[
            f(\mathbf{x}-\mathbf 1+\mathbf e_i)-f(\mathbf{x}-\mathbf 1)
        \bigr]\\
        &\ge
        f(\mathbf{x}+\mathbf e_i)-f(\mathbf{x}),
    \end{aligned}
    \]
    where the last two inequalities use the induction hypothesis and Lemma~\ref{lemma:convex_supermodular_non_increasing}.
\end{itemize}

Thus \eqref{eq:diagonal_bound} holds for $k+1$, completing the induction.
Since $B_{k+1}$ is non-decreasing for every finite horizon, if matching is weakly optimal at $\mathbf{x}\in\mathcal F$, then $B_{k+1}(\mathbf{x})\ge0$ and therefore $B_{k+1}(\mathbf{y})\ge0$ for every $\mathbf{y}\in\mathcal F$ with $\mathbf{y}\succeq\mathbf{x}$.
Hence matching is weakly optimal at every such $\mathbf{y}$, and the match-on-tie selector is monotone.
The argument is finite-horizon.
Discounted or truncated infinite-horizon variants can inherit the same monotonicity under the usual contraction or limiting conditions.
\hfill\Halmos

\subsection{Proof of Proposition~\ref{prop:fluid_balance}}
\label{app:proof_illustrative}

In a stationary fluid operating point, the matching rate must equal the arrival rate of each type, so $\bar\mu^*=\lambda/2$.
For positive fluid queue levels $(\bar X_1,\bar X_2)$, the steady-state cost rate is
\begin{equation}
\label{eq:steady_state_fluid_cost}
\bar J(\bar X_1,\bar X_2)
=
c(\bar X_1+\bar X_2)
+
\frac{\lambda}{2}
\frac{\kappa}{(\bar X_1\bar X_2)^\beta}.
\end{equation}
For a fixed average queue level $\bar X=(\bar X_1+\bar X_2)/2$, the arithmetic-geometric mean inequality gives
\[
\bar X_1\bar X_2\le \bar X^2,
\]
with equality if and only if $\bar X_1=\bar X_2=\bar X$.
Because $\beta>0$, replacing any asymmetric pair $(\bar X_1,\bar X_2)$ by $(\bar X,\bar X)$ weakly reduces the matching-cost term and leaves the waiting-cost term unchanged.
Thus the unique minimizer must be symmetric.

It remains to minimize the one-dimensional cost rate
\begin{equation}
\label{eq:symmetric_steady_state_fluid_cost}
\bar J(\bar X)
=
2c\bar X
+
\frac{\lambda}{2}\frac{\kappa}{\bar X^{2\beta}},
\qquad
\bar X>0.
\end{equation}
This objective diverges as $\bar X\downarrow0$ and as $\bar X\to\infty$, so any stationary point is a global minimizer if it is unique.
Differentiating,
\[
\bar J'(\bar X)
=
2c-\beta\lambda\kappa\bar X^{-(2\beta+1)}.
\]
The first-order condition gives
\[
(\bar X^*)^{2\beta+1}
=
\frac{\beta\lambda\kappa}{2c},
\qquad
\bar X^*
=
\left(\frac{\beta\lambda\kappa}{2c}\right)^{\frac{1}{2\beta+1}}.
\]
Moreover,
\[
\bar J''(\bar X)
=
\beta\lambda\kappa(2\beta+1)\bar X^{-(2\beta+2)}
>
0,
\]
so $\bar X^*$ is the unique minimizer.
At this point,
\[
W_{\mathrm{rate}}^*=2c\bar X^*,
\qquad
M_{\mathrm{rate}}^*=\frac{\lambda}{2}\frac{\kappa}{(\bar X^*)^{2\beta}}.
\]
Using $(\bar X^*)^{2\beta+1}=\beta\lambda\kappa/(2c)$,
\[
\frac{W_{\mathrm{rate}}^*}{M_{\mathrm{rate}}^*}
=
\frac{2c\bar X^*}{(\lambda/2)\kappa(\bar X^*)^{-2\beta}}
=
\frac{4c(\bar X^*)^{2\beta+1}}{\lambda\kappa}
=
2\beta.
\]

\hfill\Halmos

\subsection{Proof of Proposition~\ref{prop:greedy-competitive-ratio}}

\proof{Proof of Proposition~\ref{prop:greedy-competitive-ratio}.}
    We construct an adversarial arrival sequence that drives the competitive ratio to infinity.
    Consider $n$ complete tuples (one agent of each type) arriving at times $t_k = k\epsilon$ for $k = 0, 1, \ldots, n-1$, where $\epsilon > 0$ is arbitrarily small.
    Let the matching-cost function be
    \[
        f(\mathbf{x})=\frac{\bar c}{\min_i x_i}, \qquad \mathbf{x}\in\mathcal F,
    \]
    for some constant $\bar c>0$.
    This function is positive, finite, and componentwise non-increasing on $\mathcal F$, so it satisfies Assumption~\ref{ass:key}; in particular, $f(k\cdot\mathbf 1)=\bar c/k$.
    
    Under the Greedy policy, at each arrival time $t_k$, exactly one agent of each type is present in the queue, since the previous tuple was matched immediately upon its arrival.
    Therefore, Greedy executes each match at state $\mathbf{1}$, incurring cost $f(\mathbf{1}) = \bar c$ per match.
    The total cost under Greedy is
    $J^{Greedy} = n \cdot \bar c + W^{Greedy}$,
    where $W^{Greedy} \leq n \cdot (\sum_i c_i) \cdot \epsilon \to 0$ as $\epsilon \to 0$.
    
    In contrast, the optimal policy waits until time $t_{n-1} = (n-1)\epsilon$ when all $n$ tuples have arrived, then executes $n$ consecutive matches.
    The $j$-th match ($j = 1, \ldots, n$) occurs at state $(n-j+1) \cdot \mathbf{1}$, incurring cost $f((n-j+1) \cdot \mathbf{1}) = \bar c/(n-j+1)$.
    The total matching cost is $\sum_{k=1}^{n} \bar c/k = \bar c \cdot H_n$, where $H_n = \sum_{k=1}^{n} 1/k$ is the $n$-th harmonic number.
    The waiting cost satisfies $W^{OPT} \leq n \cdot (\sum_i c_i) \cdot (n-1)\epsilon \to 0$ as $\epsilon \to 0$.
    
    Comparing these costs, as $\epsilon \to 0$, the competitive ratio satisfies
    \begin{align*}
        \frac{J^{Greedy}}{J^{OPT}} \to \frac{n \cdot \bar c}{\bar c \cdot H_n} = \frac{n}{H_n}.
    \end{align*}
    Since $H_n \sim \ln n$ as $n \to \infty$, we have $n / H_n \to \infty$.
    Therefore, the competitive ratio of the Greedy policy is unbounded.
\hfill\Halmos
\endproof

\subsection{Proof of Proposition~\ref{prop:threshold-competitive-ratio}}

\proof{Proof of Proposition~\ref{prop:threshold-competitive-ratio}.}
    We construct an adversarial arrival sequence that drives the competitive ratio to infinity.
    Consider a threshold policy with parameter $\theta \geq 2$, which triggers a match only when $\min_i X_i \geq \theta$.
    
    At time $t = 0$, suppose $\theta - 1$ complete tuples arrive, so the system state is $(\theta - 1) \cdot \mathbf{1}$.
    Since $\min_i X_i = \theta - 1 < \theta$, the threshold policy does not match and waits for additional arrivals.
    At time $t = T$ for some $T > 0$, one additional complete tuple arrives.
    The resulting finite arrival sequence is balanced, with $\theta$ arrivals of each type.
    
    Under the Threshold policy, all $(\theta - 1) \cdot N$ agents wait in the system for the entire interval $[0, T]$.
    The waiting cost incurred is $W^{Threshold} = (\theta - 1) \cdot (\sum_{i=1}^{N} c_i) \cdot T$.
    Regardless of how the policy clears agents after time $T$, its total cost satisfies
    $J^{Threshold} \geq (\theta - 1) \cdot (\sum_{i=1}^{N} c_i) \cdot T$.
    
    In contrast, the optimal policy matches the first $\theta - 1$ tuples immediately at time $t = 0$ and matches the final tuple when it arrives at time $T$.
    Its total cost satisfies
    \[
        J^{OPT} \leq \sum_{k=1}^{\theta-1} f(k \cdot \mathbf{1}) + f(\mathbf{1}),
    \]
    which is finite and independent of $T$.
    
    The competitive ratio satisfies
    \begin{align*}
        \frac{J^{Threshold}}{J^{OPT}} \geq \frac{(\theta - 1) \cdot (\sum_{i=1}^{N} c_i) \cdot T}{\sum_{k=1}^{\theta-1} f(k \cdot \mathbf{1})+f(\mathbf{1})} \to \infty \quad \text{as } T \to \infty.
    \end{align*}
    Therefore, the competitive ratio of the Threshold policy is unbounded.
\hfill\Halmos
\endproof

\subsection{Proof of Proposition~\ref{prop:lower-bound}}

\proof{Proof of Proposition~\ref{prop:lower-bound}.}
    It is enough to consider balanced complete-tuple arrivals.
    Rescale time so that the waiting-cost rate of one complete tuple is one, i.e.,
    \[
        w(\mathbf{1})=\sum_{i=1}^N c_i=1.
    \]
    For a feasible state $\mathbf{x}\in\mathcal F$, let
    \[
        q(\mathbf{x})=\min_i x_i
    \]
    denote the number of complete tuples present in the pool.
    Fix $0<\delta\leq 1$ and an integer $m\geq1$.
    Define
    \[
        f_\delta(\mathbf{x})
        =
        \begin{cases}
        1, & q(\mathbf{x})=1,\\
        \delta, & q(\mathbf{x})\geq 2.
        \end{cases}
    \]
    This matching-cost function is strictly positive and componentwise non-increasing on $\mathcal F$.
    Set $f=f_\delta$ for the constructed instance.

    Fix a deterministic online algorithm, denoted ALG.
    The adversary releases complete tuples one at a time.
    At time $a_1=0$, one complete tuple arrives.
    After the $j$-th tuple arrives at time $a_j$, suppose no future arrivals occur, and let $d_j$ be the time ALG would take, measured from $a_j$, to empty the system.
    If ALG would never empty the system, set $d_j=\infty$.

    The adversary uses the following adaptive rule.
    If $d_j\geq 1$, the adversary stops after the $j$-th tuple.
    If $d_j<1$, the adversary waits until ALG empties the system and then releases the next complete tuple at time
    \[
        a_{j+1}=a_j+d_j+\epsilon,
    \]
    where $\epsilon>0$ is arbitrarily small.
    The adversary continues this process until either the stopping condition occurs or $m$ tuples have been released.

    Suppose first that the stopping condition does not occur before $m$ tuples have been released.
    Along the realized path, each new tuple arrives only after ALG has emptied the system.
    Hence ALG always matches a singleton complete tuple, paying matching cost $f_\delta(\mathbf{1})=1$ each time, and incurs waiting cost $d_j$ on tuple $j$.
    Therefore,
    \[
        J^{ALG}\geq
        L_m
        :=
        m+\sum_{j=1}^m d_j .
    \]

    Now consider OPT on the same arrival sequence.
    OPT can batch adjacent complete tuples using one of the two pairings
    \[
        (1,2),(3,4),(5,6),\ldots
    \]
    or
    \[
        (2,3),(4,5),(6,7),\ldots,
    \]
    with at most two singleton tuples left unmatched by pairs.
    If OPT batches tuple $j$ with tuple $j+1$, tuple $j$ waits for $d_j+\epsilon$ time units, and the two tuples are cleared at matching cost $1+\delta$: the first match is made from a state with two complete tuples at cost $\delta$, and the second from a singleton state at cost $1$.
    Averaging over the two pairings, one of them has cost at most
    \[
        \frac{m+1}{2}
        +
        \frac{1}{2}\sum_{j=1}^{m-1}d_j
        +
        \frac{\delta m}{2}
        +
        O(m\epsilon).
    \]
    Letting $\epsilon\downarrow0$ and using $L_m=m+\sum_{j=1}^m d_j$, we obtain
    \[
        J^{OPT}
        \leq
        \frac{L_m+1}{2}
        +
        \frac{\delta L_m}{2},
    \]
    where we used $L_m\geq m$.
    Thus
    \[
        \frac{J^{ALG}}{J^{OPT}}
        \geq
        \frac{L_m}{(L_m+1)/2+\delta L_m/2}
        =
        \frac{2}{1+\delta+1/L_m}
        \geq
        \frac{2}{1+\delta+1/m}.
    \]

    It remains to handle the case in which the adversary stops because $d_j\geq 1$ for some $j\leq m$.
    If $d_j=\infty$, then ALG has infinite cost and the lower bound is immediate.
    Thus suppose $1\leq d_j<\infty$.
    Applying the same adjacent-pairing argument to the prefix of $j$ tuples gives
    \[
        J^{OPT}
        \leq
        \frac{j+1}{2}
        +
        \frac{1}{2}\sum_{r=1}^{j-1}d_r
        +
        \frac{\delta j}{2}.
    \]
    Since $d_j\geq1$ and
    \[
        L_j=j+\sum_{r=1}^j d_r,
    \]
    the possible extra singleton in the pairing bound is offset by ALG's final delay:
    \[
        \frac{j+1}{2}
        +
        \frac{1}{2}\sum_{r=1}^{j-1}d_r
        +
        \frac{\delta j}{2}
        \leq
        \frac{1+\delta}{2}L_j.
    \]
    Hence
    \[
        J^{OPT}
        \leq
        \frac{1+\delta}{2}L_j.
    \]
    Because ALG pays at least $L_j$ on this prefix,
    \[
        \frac{J^{ALG}}{J^{OPT}}
        \geq
        \frac{2}{1+\delta}.
    \]

    Finally, fix $\eta>0$.
    If $\eta\geq2$, the claim is immediate, so assume $0<\eta<2$.
    Choose $0<\delta\leq1$ and an integer $m$ large enough that
    \[
        \delta+\frac{1}{m}\leq \frac{\eta}{2-\eta}.
    \]
    In the no-early-stop case, the ratio is at least
    \[
        \frac{2}{1+\delta+1/m}\geq 2-\eta.
    \]
    In the early-stop case, the ratio is at least $2/(1+\delta)$, which is also at least $2-\eta$ by the same choice of $\delta$.
    Therefore, for every deterministic online algorithm, the adversary constructs a finite balanced complete-tuple arrival sequence with competitive ratio at least $2-\eta$.
\hfill\Halmos
\endproof

\subsection{Proof of Proposition~\ref{prop:fixed-consumption}}

\proof{Proof of Proposition~\ref{prop:fixed-consumption}.}
    \textit{Part 1: Upper bound.}
    The proof follows the proof of Theorem~\ref{thm:main} with $\mathbf 1$ replaced by $\mathbf r$.
    We give the argument to make the required invariant explicit.
    Let $\tau_k$ be the time of CB's $k$-th match, let $s_k$ be the time of OPT's $k$-th match, and let $W_k$ be CB's accumulated waiting cost over $(\tau_{k-1},\tau_k]$.
    Let $\mathbf x_k$ and $\mathbf u_k$ be the pre-match states for CB's and OPT's $k$-th matches, respectively, and write $M_k=f_{\mathbf r}(\mathbf x_k)$.
    The CB trigger gives
    \[
        W_k+M_k\le (1+\alpha)W_k.
    \]
    
    If $s_k\ge \tau_k$, then throughout CB's $k$-th intermatch interval CB has completed $k-1$ matches and OPT has completed at most $k-1$ matches.
    Hence, outside event times of measure zero,
    \[
        \mathbf X^{CB_\alpha(\mathbf r)}(t)
        =
        \mathbf A(t)-(k-1)\mathbf r,
        \qquad
        \mathbf X^{OPT(\mathbf r)}(t)
        =
        \mathbf A(t)-m^{OPT}(t)\mathbf r,
        \quad m^{OPT}(t)\le k-1.
    \]
    Because $\mathbf r\succeq\mathbf 0$, OPT's state componentwise dominates CB's state, so OPT's waiting-cost rate is at least CB's.
    Therefore CB's $k$-th cost is at most
    \[
        (1+\alpha)\int_{\tau_{k-1}}^{\tau_k}w(\mathbf X^{OPT(\mathbf r)}(t))\,dt.
    \]
    
    If $s_k<\tau_k$, then both policies have completed $k-1$ matches immediately before their respective $k$-th matches, and the cumulative arrival vector at $s_k$ is componentwise no larger than at $\tau_k$.
    Thus $\mathbf u_k\preceq\mathbf x_k$, and monotonicity gives $M_k\le f_{\mathbf r}(\mathbf u_k)$.
    It remains to bound $W_k$.
    If CB's $k$-th match is caused by a timer, then $M_k=\alpha W_k$, so $W_k\le f_{\mathbf r}(\mathbf u_k)/\alpha$.
    If it is caused by an arrival, let $\mathbf y_k$ be CB's state immediately before that triggering arrival.
    We claim that $\mathbf y_k\in\mathcal F_{\mathbf r}$.
    Otherwise, for some type $i$, $y_{k,i}<r_i$.
    Since CB has already consumed $(k-1)r_i$ type-$i$ agents before its $k$-th match, this implies $A_i(\tau_k^-)<k r_i$.
    But OPT executed its $k$-th match earlier, so feasibility of that match requires $A_i(s_k)\ge k r_i$, contradicting $s_k<\tau_k$.
    Hence $\mathbf y_k\in\mathcal F_{\mathbf r}$.
    Timer-priority tie-breaking implies that CB's trigger was false immediately before the triggering arrival:
    \[
        f_{\mathbf r}(\mathbf y_k)>\alpha W_k.
    \]
    Because $\mathbf u_k\preceq\mathbf y_k$, monotonicity gives
    \[
        f_{\mathbf r}(\mathbf u_k)\ge f_{\mathbf r}(\mathbf y_k)>\alpha W_k.
    \]
    Hence CB's $k$-th cost is at most $(1+1/\alpha)f_{\mathbf r}(\mathbf u_k)$.
    
    Summing over the two cases completes the charging argument.
    The CB intermatch intervals charged to OPT's waiting cost are disjoint, and the early-OPT cases charge distinct OPT matching costs.
    Hence the total charge is at most OPT's total cost, exactly as in Theorem~\ref{thm:main}.
    This proves the stated competitive ratio.

    \textit{Part 2: Tightness.}
    Fix any $\eta>0$ and any deterministic online policy $\pi$.
    The proof reduces the fixed-consumption problem to the lower-bound construction in Proposition~\ref{prop:lower-bound} by treating one full vector $\mathbf r$ as a single block.
    Define
    \[
        q_{\mathbf r}(\mathbf x)
        :=
        \min_i\left\lfloor \frac{x_i}{r_i}\right\rfloor,
    \]
    the number of complete blocks available in state $\mathbf x$.
    Rescale time so that the waiting-cost rate of one block is one:
    \[
        w(\mathbf r)=\sum_i c_i r_i=1.
    \]
    For $0<\delta\le 1$, set
    \[
        f_{\mathbf r,\delta}(\mathbf x)
        =
        \begin{cases}
        1, & q_{\mathbf r}(\mathbf x)=1,\\
        \delta, & q_{\mathbf r}(\mathbf x)\ge 2.
        \end{cases}
    \]
    This cost function is positive and componentwise non-increasing on $\mathcal F_{\mathbf r}$.
    The adversary releases one block at a time.
    After the $j$-th block arrives, suppose no future arrivals occur, and let $d_j$ be the time the deterministic policy would take, measured from that arrival time, to empty the system.
    If $d_j\ge1$, the adversary stops; if $d_j<1$, the adversary waits until the policy empties the system and then releases the next block after an arbitrarily small delay $\epsilon$.
    If the stopping condition never occurs before $m$ blocks have been released, the online policy clears singleton blocks along the realized path and pays at least
    \[
        L_m=m+\sum_{j=1}^m d_j .
    \]
    As in Proposition~\ref{prop:lower-bound}, OPT can use one of the two adjacent pairings of released blocks.
    Pairing block $j$ with block $j+1$ makes block $j$ wait for $d_j+\epsilon$ and then clears the two-block state at cost $\delta$ followed by the remaining singleton block at cost $1$.
    Averaging over the two adjacent pairings gives the same bound
    \[
        J^{OPT(\mathbf r)}
        \le
        \frac{L_m+1}{2}
        +
        \frac{\delta L_m}{2}
        +
        O(m\epsilon).
    \]
    Letting $\epsilon\downarrow0$ yields
    \[
        \frac{J^\pi}{J^{OPT(\mathbf r)}}
        \ge
        \frac{2}{1+\delta+1/m}
    \]
    in the no-early-stop case.
    If the adversary stops at some $j$ with $d_j\ge1$, the same prefix argument as in Proposition~\ref{prop:lower-bound} gives
    \[
        \frac{J^\pi}{J^{OPT(\mathbf r)}}\ge \frac{2}{1+\delta}.
    \]
    Choosing $\delta$ small and $m$ large enough gives a ratio at least $2-\eta$.
\hfill\Halmos
\endproof

\subsection{Proof of Proposition~\ref{prop:nphard}}
\label{app:proof_nphard}

\proof{Proof.}
Consider a static instance where all agents arrive at $t = 0$, waiting costs are zero ($c_i = 0$ for all $i$), and the controller must partition agents into disjoint $N$-tuples to minimize total matching cost.
For $N = 3$, this reduces to the 3-Dimensional Matching (3DM) problem.

Let there be $n$ agents of each type with attribute sets $\mathcal{Q}_1 = \{q_{1,1}, \ldots, q_{1,n}\}$, $\mathcal{Q}_2 = \{q_{2,1}, \ldots, q_{2,n}\}$, $\mathcal{Q}_3 = \{q_{3,1}, \ldots, q_{3,n}\}$.
Let $T \subseteq \mathcal{Q}_1 \times \mathcal{Q}_2 \times \mathcal{Q}_3$ be a set of feasible triples.
Define the matching cost as
$$
f(q_{1,i}, q_{2,j}, q_{3,k}) = 
\begin{cases}
0 & \text{if } (q_{1,i}, q_{2,j}, q_{3,k}) \in T, \\
1 & \text{otherwise}.
\end{cases}
$$
Finding a matching with total cost zero is equivalent to finding a perfect 3-dimensional matching, which is NP-complete \citep{Karp1972}.
Since the dynamic problem generalizes this static instance, computing optimal policies in the attribute-based model is NP-hard.
\hfill\Halmos
\endproof

\subsection{Proof of Proposition~\ref{prop:adversarial attribute}}
\label{app:proof_adversarial_attribute}

\proof{Proof of Proposition~\ref{prop:adversarial attribute}.}
    We construct an adversarial strategy that forces any deterministic online algorithm to have unbounded competitive ratio by exploiting attribute-dependent costs.
    Consider a two-sided matching problem ($N = 2$) in the general attribute-based model of Section~\ref{sec:attribute-model}, where the matching cost for a pair of attributes $(q_1, q_2) \in \mathcal{Q} \times \mathcal{Q}$ is given by a function $f(q_1, q_2)$, and the waiting cost rate is $c > 0$ per agent per unit time.
    
    Fix four attributes $a_1, a_2, b_1, b_2 \in \mathcal{Q}$ such that
    \begin{align*}
        f(a_1, b_1) &= M_1, \\
        f(a_1, b_2) &= \delta, \\
        f(a_2, b_1) &= \delta, \\
        f(a_2, b_2) &= M_2,
    \end{align*}
    where $M_1 > 0$, $0<\delta<M_1$, and $M_2 \gg M_1$ will be chosen by the adversary.
    It is straightforward to verify that this definition of $f(\cdot, \cdot)$ satisfies Assumption~\ref{ass:general_attribute-based_model}.
    
    At time $t = 0$, a type-1 agent with attribute $a_1$ and a type-2 agent with attribute $b_1$ arrive.
    Let ALG be any deterministic online algorithm, and let $\tau \geq 0$ be the time at which ALG decides to match this pair.
    If ALG never matches them, then Case 2 below gives infinite waiting cost, so suppose $\tau<\infty$.
    At time $\tau + \epsilon$ (for arbitrarily small $\epsilon > 0$), the adversary adapts based on ALG's choice of $\tau$ and decides whether to send an additional pair of agents.
    
    Case 1 (adversary sends new agents): At time $\tau + \epsilon$, the adversary sends a type-1 agent with attribute $a_2$ and a type-2 agent with attribute $b_2$.
    ALG has already matched the initial pair at time $\tau$, incurring matching cost $M_1$ and waiting cost $2c\tau$, for a total of $M_1 + 2c\tau$.
    When the new pair $(a_2, b_2)$ arrives, ALG must eventually match them as well, incurring an additional cost of at least $M_2$.
    Thus,
    \[
        J^{ALG} \geq M_1 + M_2 + 2c\tau.
    \]
    In contrast, the optimal offline policy waits until time $\tau + \epsilon$, when all four agents are present, and matches $(a_1, b_2)$ and $(a_2, b_1)$.
    The offline policy pays the waiting cost of the initial pair plus the two positive cross-pair matching costs:
    \[
        J^{OPT} = 2c(\tau + \epsilon) + 2\delta \to 2c\tau + 2\delta \quad \text{as } \epsilon \to 0.
    \]
    Hence, since $\epsilon$ can be chosen arbitrarily small, the limiting competitive ratio in this case satisfies
    \[
        \lim_{\epsilon\downarrow 0}\rho_1(\tau) \;\geq\; \frac{M_1 + M_2 + 2c\tau}{2c\tau+2\delta}.
    \]
    
    Case 2 (adversary sends no further agents): 
    No additional agents arrive after $t = 0$.
    ALG matches the initial pair at time $\tau$, incurring total cost
    \[
        J^{ALG} = M_1 + 2c\tau.
    \]
    The optimal offline policy matches them immediately at $t = 0$, incurring cost
    \[
        J^{OPT} = M_1.
    \]
    The competitive ratio in this case is
    \[
        \rho_2(\tau) = \frac{M_1 + 2c\tau}{M_1}
        = 1 + \frac{2c\tau}{M_1}.
    \]
    
    Combining the two cases, and taking $\epsilon$ arbitrarily small in Case 1, it suffices to lower-bound
    \[
        \max\left\{\frac{M_1+M_2+2c\tau}{2c\tau+2\delta},\; 1 + \frac{2c\tau}{M_1}\right\}.
    \]
    To see that this lower bound is unbounded, fix any target ratio $R\geq 1$ and let $r=M_2/M_1$.
    Choose $r>R(R-1+2\delta/M_1)$.
    If $1+2c\tau/M_1>R$, then Case 2 already gives ratio larger than $R$.
    Otherwise, $2c\tau/M_1\le R-1$, and the first term in the bound above satisfies
    \[
        \frac{M_1+M_2+2c\tau}{2c\tau+2\delta}
        \geq
        \frac{M_2}{2c\tau+2\delta}
        \geq
        \frac{r}{R-1+2\delta/M_1}
        >
        R.
    \]
    Since the inequality is strict, the adversary can choose $\epsilon>0$ small enough in Case 1 to preserve the ratio above $R$.
    Therefore, no deterministic online algorithm achieves a bounded competitive ratio when the adversary can control both arrival times and attributes.
\hfill\Halmos
\endproof

\end{document}